\documentclass{patera}

\usepackage[cp1251]{inputenc}
\include{diagxy}
\input{xy}
\xyoption{all}\CompileMatrices

\def\C{\mathbb C}
\def\R{\mathbb R}
\def\Z{\mathbb Z}
\def\N{\mathbb N}

\def\r{\rangle}
\def\l{\langle}

\def\w{\omega}

\newcommand{\cw}{{\check\omega}}

\newcommand{\comment}[1]{}

\newtheorem*{lemma*}{Lemma}

\newtheorem*{conjecture*}{Conjecture}

{\theoremstyle{definition}
\newtheorem{definition}{Definition}
\newtheorem{example}{Example}

\newtheorem*{note*}{Note}
}

\allowdisplaybreaks

\begin{document}
\begin{flushleft}
\Large \bf Orthogonal polynomials of compact simple Lie groups:\\
Branching rules for polynomials
\end{flushleft}
\medskip

\begin{flushleft}
M. Nesterenko$^\S$,
J. Patera$^\dag$,
M. Szajewska$^\ddag$, and A. Tereszkiewicz$^\ddag$
\medskip

\end{flushleft}
\noindent $^\S$~Institute of Mathematics of NAS of Ukraine, 3 Tereshchenkivs'ka Str., Kyiv-4, 01601, Ukraine;\\
\noindent $^\dag$~CRM, Universit\'e de Montr\'eal, C.P.6128-Centre ville, Montr\'eal, Canada;\\
\noindent $^\ddag$~Institute of Mathematics, University of Bialystok, Akademicka 2, PL-15-267, Bialystok, Poland.
\medskip

\noindent~E-mail:            $^\S$maryna@imath.kiev.ua,  $^\dag$patera@crm.umontreal.ca,\\
\phantom{\noindent~E-mail:\ }$^\ddag$marzena@math.uwb.edu.pl,\; a.tereszkiewicz@uwb.edu.pl

{\vspace{9mm}\par\noindent\hspace*{8mm}\parbox{140mm}
{\small
Polynomials in this paper are defined starting from a compact semisimple Lie group.
A~known classification of maximal, semisimple subgroups of simple Lie groups is used to select the cases to be considered here.
A~general method is presented and all the cases of rank $\leq3$ are explicitly studied.
We derive the polynomials of simple Lie groups $B_3$ and $C_3$ as they are not available elsewhere.
The results point to far reaching Lie theoretical connections to the theory of multivariable orthogonal polynomials.
\par\vspace{7mm}}}

\section{Introduction}
The main purpose of the paper is to demonstrate, describe and illustrate homomorphic relations  (also called reduction or branching) between families of polynomials with different underlying Lie groups.
The polynomials can be viewed as multivariable generalizations of classical Chebyshev polynomials of one variable
or as subfamilies of multivariable Macdonald polynomials~\cite{Macdonald1}.
The systematic study of such relations became possible after several families of polynomials in $n$ variables were constructed~\cite{NPT} via semisimple Lie groups of rank~$n$.

The relations studied here are consequences of maximal inclusions of semisimple Lie groups in simple compact Lie groups. They parallel familiar branching rules for finite dimensional representations of corresponding Lie algebras (see for example \cite{McP}), but cannot be obtained from them in any direct way.  The technique exploited in the computation of the branching rules for polynomials is the adaptation of the method previously used  in \cite{LNP} and \cite{McP}. Related problems are the computation of branching rules for orbit functions \cite{KP1,KP2} and Weyl group orbits~\cite{LNP}.

The problems share certain tools, namely the projection matrices for weights.
All the cases of interest to us were classified half a century ago by E.B.~Dynkin.
In this paper we  describe the general principle of the method, and consider all the specific cases with rank $n\leqslant3$.

In the literature~\cite{NPT} one finds sufficiently many explicit examples of polynomials for the groups of rank $n=2$, but only those of $A_3$ for $n=3$.
Therefore we start by computing the polynomials of the group $B_3$ and $C_3$.
As in \cite{NPT}, we take the orbit functions of the three fundamental weights  for either of the two Lie groups as our polynomial variables. Such a substitution imposes transformations of the orthogonality domains. For the orbit functions of $B_3$ and $C_3$ these were tetrahedra inscribed in a cube. The polynomial substitution transforms them into the domains $\mathfrak F$ shown in Figures~\ref{fig_b3} and \ref{fig_c3}.

New here is description of the continuous and discrete orthogonality of the polynomials within domains $\mathfrak F$.

We start from the $C$-functions of \cite{KP1} and from the $S$-functions of \cite{KP2} (see also \cite{P,MP}),
by specializing them to three variables and converting them into polynomials.  The substitution of variables providing the conversion was introduced in \cite{NPT}.
The simplest 1-variable version of such a conversion is found in the context of Chebyshev polynomials~\cite{R}  and their generalization to two variables \cite{NPT}, \cite{Koor}.
The background to our work is at the crossroad of the theory of compact semisimple Lie groups particularly the properties of their characters, the theory of special functions of mathematical physics, and the orthogonal polynomials of many variables.
From the Lie theory we take the uniformity of description of properties of each simple Lie group and of their characters.
The price to pay for that is the need to work with non-orthogonal bases and with character functions of ever increasing complexity the larger the weights get.
By~working with the orbit functions instead, we circumvent the problem of characters while still working with their $W$-invariant constituents.
Contact with the theory of special functions is made through properties of orbit functions that are symmetric and skew-symmetric with respect to boundaries of $\mathfrak F$, called the $C$- and $S$-functions respectively \cite{KP1,KP2}.
The two families of functions have been part of Lie theory for almost a century.
Their properties as special functions, particularly their discrete orthogonality, were recognized only a few decades ago \cite{MP87}.
Although the characters could also be viewed as special functions of mathematical physics, their complexity disqualifies them from almost all applications.
To the best of our knowledge the idea to see the root systems as the backbone of the theory of orthogonal polynomials of many variables comes from~\cite{Macdonald}.
Here we use the theory of simple Lie groups in order to construct and reduce multivariate orthogonal polynomials.
In retrospect, the results of \cite{Koor} are based on group $A_2$, the results of \cite{LX} on~$A_n$.
The approach exploits complete reducibility of products of orbit functions in order to build polynomials from the lowest few.
The orbit functions of fundamental weights of the Lie group become the polynomial variables. Let us emphasize that there are alternatives to our approach that are untested so far. Products of characters are also completely decomposable into their sums. Therefore a similar recursive procedure would build the polynomials as functions of variables that are characters of fundamental representations of the underlying Lie group. Due to the basic role of the characters, this version of our method can be preferable for some problems. However, the complex structure of the characters, as opposed to orbit functions, makes it practically more cumbersome.
The~two approaches coincide for simple Lie groups $A_n$. In such cases the characters of the fundamental weights are equal to the orbit functions of the same weights.

\section{Preliminaries}
The notions of polynomials under consideration in $n$ variables depend essentially on the underlying semisimple Lie group~$G$ of rank $n$.
This section is intended to fix notation and terminology and to recall the definitions and some properties of orbit functions.
Additional information on this subject can be found, for example, in~\cite{KassMoodyPateraSlansky1990, HP,KP1,KP2}.

\subsection{Bases, Weyl group and orbit functions}
Let $\R^n$ be the real Euclidean space spanned by the simple roots of a~simple Lie group~$G$ (equivalently, Lie algebra).
The basis of the simple roots and the basis of fundamental weights are hereafter referred to as the $\alpha$-basis and $\omega$-basis, respectively.
The two bases are linked by the Cartan matrix $\mathfrak{C}$ in the following way $\alpha=\mathfrak{C}\omega$ where
\begin{gather*}
\mathfrak{C}:=(\mathfrak{C}_{jk})= \left(\frac{2\langle\alpha_j,\alpha_k\rangle}{\langle\alpha_k,\alpha_k\rangle}\right),
\qquad
\text{hereafter} \quad j,k=1,2,\dots,n.
\end{gather*}
The Cartan matrix provides, in principle, all the information needed about $G$.
The same data about group~$G$ can be taken from the Coxeter-Dynkin diagrams, see e.g.~\cite{McP}.

We also use the convention that for the long roots $\alpha_k$ the inner product $\langle\alpha_k ,\alpha_k\rangle=2$,
and we introduce bases dual to $\alpha$- and $\omega$-bases, denoted here as $\check{\alpha}$- and $\check{\omega}$-bases, respectively.
The dual bases are fixed by the relations
\begin{gather*}\label{duals}
\check{\alpha}_j=\frac{2\alpha_j}{\langle\alpha_j ,\alpha_j\rangle},
\qquad
\check{\omega}_j=\frac{2\omega_j}{\langle\alpha_j ,\alpha_j\rangle},
\qquad
\langle\alpha_j ,\check{\omega}_k\rangle=\langle\check{\alpha}_j ,\omega_k\rangle=\delta_{jk},
\end{gather*}
Occasionally it is also useful to work with the orthonormal basis $\{e_1,e_2,\ldots, e_n\}$ of $\R^n$.

Now we can form the root lattice $Q$ and the weight lattice $P$  of $G$
by all integer linear combinations of the $\alpha$-basis and
$\omega$-basis,
\begin{gather*}
Q=\Z\alpha_1+\Z\alpha_2+\cdots+\Z\alpha_n,\qquad
P=\Z\omega_1+\Z\omega_2+\cdots+\Z\omega_n.
\end{gather*}
In the weight lattice $P$, we define the cone of dominant weights
$P^+$ and its subset of strictly dominant weights $P^{++}$
\begin{gather*}
P\;\supset\; P^+=\Z^{\ge 0}\omega_1+\cdots+\Z^{\ge 0}\omega_n
\;\supset\; P^{++}=\Z^{>0}\omega_1+\cdots+\Z^{>0}\omega_n.
\end{gather*}
Analogously dual latices $\check{Q}$ and $\check{P}$ are defined as follows
\begin{gather*}
\check{Q}=\Z\check{\alpha}_1+\Z\check{\alpha}_2+\cdots+\Z\check{\alpha}_n,\qquad
\check{P}=\Z\check{\omega}_1+\Z\check{\omega}_2+\cdots+\Z\check{\omega}_n.
\end{gather*}

Weyl group $W(G)$ is the finite group generated by reflections in $(n-1)$-dimensional hyperplanes orthogonal to simple roots, having the origin as their common point, and referred to as elementary reflections $r_{\alpha_j}=r_j$, $j=1,\dots, n$.

The orbit of $W$ containing the (dominant) point $\lambda\in P^+\subset \R^n$ is written as  $W_\lambda$.
The size of $W_\lambda$ is denoted by  $|W_\lambda|$ (it is the number of points in $W_\lambda$), and order of the Weyl group is denoted by $|W|$.

The fundamental region $F(G)\subset \R^n$ is the convex hull of the vertices
$\{0,\frac{\check{\omega}_1}{m_1},\ldots, \frac{\check{\omega}_n}{m_n}\}$, where $m_j$, $j=\overline{1,n}$ are marks of the highest root $\xi$ of the root system.

In this paper we mainly deal with the simple Lie groups of rank three,
and in the Appendix we present all necessary information about such groups,
i.e., their Cartan matrices, Weyl group orbits, highest roots, Weyl orbit sizes and fundamental regions.
The above brings us to definitions of symmetric and antisymmetric orbit functions.

\begin{definition}\label{C}
The $C$-function $C_{\lambda}(x)$ of $G$ is defined as
\begin{gather*}\label{def_c-function1}
C_\lambda(x) := \sum_{\mu\in W_\lambda(G)} e^{2\pi i \l\mu, x\r},
\qquad x\in\R^n,\quad \lambda\in P^+.
\end{gather*}
\end{definition}

\begin{definition}
The $S$-function $S_{\lambda}(x)$ is defined as
\begin{gather*}\label{def_s-function1}
S_\lambda(x) := \sum_{\mu\in W_\lambda(G)} (-1)^{p(\mu)}e^{2\pi i\l\mu,x\r},
\qquad x\in\R^n,\quad \lambda\in P^{++}\,.
\end{gather*}
where $p(\mu)$ is the number of elementary reflections  necessary to obtain $\mu$ from~$\lambda$.
\end{definition}
The same $\mu$ can be obtained by different successions of reflections,
but all shortest routes from $\lambda$ to $\mu$ will have same parity length, so $S$-functions are well defined.

In this paper, we always suppose that $\lambda,\ \mu\in P$ are given in $\w$-basis and $x\in\R^n$ is given in $\check{\alpha}$ basis, namely
$\lambda=\sum\limits^n_{j=1}\lambda_j\w_j$, \
$\mu=\sum\limits^n_{j=1}\mu_j\w_j$, $\lambda_j,\ \mu_j\in\Z$ and
$x=\sum\limits^n_{j=1}x_j\check{\alpha}_j$, \ $x_j\in \R$.
Therefore the orbit functions have the following forms
\begin{gather}
C_\lambda(x)
= \sum_{\mu\in W_\lambda} e^{2\pi i \sum\limits^n_{j=1}\mu_jx_j}
= \sum_{\mu\in W_\lambda} \prod\limits^n_{j=1} e^{2\pi i \mu_jx_j},
\label{Cdef}
\\
S_\lambda(x)
= \sum_{\mu\in W_\lambda} (-1)^{p(\mu)}e^{2\pi i \sum\limits^n_{j=1}\mu_jx_j}
= \sum_{\mu\in W_\lambda} (-1)^{p(\mu)}\prod\limits^n_{j=1} e^{2\pi i \mu_jx_j}\label{Sdef}.
\end{gather}

The introduced orbit functions have many useful properties,
e.g., continuity, orthogonality, symmetry (antisymmetry) with respect to the boundary of~$F$, eigenfunctions of the differential operators, etc.
(for details see~\cite{KP1, KP2}).

\subsection{Discretization of orbit functions}\label{sec_discretization}
Both $C$- and $S$-families of orbit functions are orthogonal and complete, which makes them perfect for the Fourier analysis.
As a lattice for the Fourier analysis we choose the refinement of the $\Z$-dual lattice to $Q$, namely $\frac 1 M \check{P}$,
where $M\in\N$.

Repeated reflections of $F(G)$ in its $(n-1)$-dimensional sides results in tiling the entire space $\R^n$ by copies of $F$,
moreover it is sufficient to consider orbit functions only on the fundamental region, therefore let us discretize $F$.

We define $F_M\subset F$, depending on an arbitrary natural number $M$ as follows
\begin{gather*}
F_M=\left\{
\frac{s_1}{M}\check\omega_1+\frac{s_2}{M}\check\omega_2+\cdots+\frac{s_n}{M}\check\omega_n\mid
s_1,\dots,s_n\in\Z^{\geq0},\ \;
\sum_{i=1}^n s_im_i\le M\in\N
\right\}.
\end{gather*}
For $S$-functions, the discretized fundamental region is the interior of $F_M$ and it has the form
\begin{gather*}
\widetilde{F}_M=\left\{
\frac{s_1}{M}\check\omega_1+\frac{s_2}{M}\check\omega_2+\cdots+\frac{s_n}{M}\check\omega_n\mid
s_1,\dots,s_n\in\Z^{>0},\ \;
\sum_{i=1}^n s_im_i\le M\in\N
\right\}.
\end{gather*}
The number of points of the grid $F_M$ (or $\widetilde{F}_M$) is denoted by $|F_M|$ (or $|\widetilde{F}_M|$).

We define the scalar product (see~\cite{HP}) of two functions $f,g:F_M(G)\rightarrow \C$ by
\begin{gather*}
\langle f, g \rangle_{F_M}=\sum\limits_{x\in {F_M}} \varepsilon(x)f(x)\overline{g(x)},
\qquad \text{where}\quad
\varepsilon(x):=|W_x|.
\end{gather*}
The same orthogonality relation holds true for $f,g:\widetilde{F}_M\rightarrow \C$.

For $C$- and $S$-functions normalized by the order of stabilizer of $\lambda$, we have:
\begin{gather}\label{Corth}
\langle C_\lambda(x), C_{\lambda'}(x) \rangle_{F_M}=
\sum\limits_{x\in {F_M}} {|W_x|} C_\lambda(x)\overline{C_{\lambda'}(x)}=
\det \mathfrak{C}\ \frac{|W|^2}{|W_\lambda|}M^n\delta_{\lambda\ \lambda'},\\
\label{Sorth}\langle S_\lambda(x), S_{\lambda'}(x) \rangle_{\widetilde{F}_M}=
\sum\limits_{x\in\widetilde{F}_M} {|W|}S_\lambda(x)\overline{S_{\lambda'}(x)}=
\det \mathfrak{C}\ |W|\ M^n\delta_{\lambda\ \lambda'}.
\end{gather}

Precise values of $|W|$, $|W_\lambda|$, $|F_M|$ and $|\widetilde{F}_M|$ for the simple Lie groups of rank three are presented in the Appendix.

\section{Orthogonal polynomials in $n$ variables}
In this section we fix notations, recall definitions of multivariate polynomials of simple Lie groups introduced in~\cite{NPT}
and explain some useful notions taken from~\cite{Dunkl} and \cite{Xu}.

The main objects of this paper are polynomials in $n$ variables
\begin{gather}
P_{k_1,\ldots,k_n}(u){=}\!\!\!\sum\limits_{j_1,\ldots,j_n=0}^{k_1,\ldots,k_n}\!a_{j_1,\ldots,j_n}u_1^{j_1}\cdots u_n^{j_n},
\quad \text{where}\quad
u:=(u_1,\ldots,u_n)\in \C^n,\;a_{j_1,\ldots,j_n}\in \R.
\end{gather}

\begin{definition}[\textit{level vector order}]\label{order}
Let $(a_1,a_2,\ldots,a_n)$ be the level vector for $G$ (see \cite{BMP})
and let
$\mathfrak{n}= (a_1,a_2,\ldots,a_n)(k_1,k_2,\dots, k_n)^t$,
and
$\mathfrak{n}'= (a_1,a_2,\ldots,a_n)(k'_1,k'_2,\dots, k'_n)^t$.

Then we say that
\begin{gather*}
u_1^{k_1}u_2^{k_2}\cdots u_n^{k_n} \succ u_1^{k_1'}u_2^{k'_2}\cdots u_n^{k_n'}
\end{gather*}
if $\mathfrak{n}>\mathfrak{n}'$ or if $\mathfrak{n}=\mathfrak{n}'$ and  the first nonzero entry in the $n$-tuple $(k_1{-}k_1',  k_2{-}k_2',\ldots,k_n{-}k_n')$ is negative.
\end{definition}

In fact the level vector order for vectors $(k_1,k_2,\ldots,k_n)$ and $(k_1',k_2',\ldots,k_n')$ coincides with the graded lexicographical order for $(a_1k_1,a_2k_2,\ldots,a_nk_n)$ and $(a_1k_1',a_2k_2',\ldots,a_nk_n')$ vectors.

As soon as the above ordering is fixed, the highest modified total degree of the monomials  $u_1^{k_1}\cdots u_n^{k_n}$
(i.e. ${\rm max}\{a_1k_1+\cdots+a_nk_n\}$)
of the polynomial $P_{k_1,\ldots,k_n}(u)$ is called the \textit{modified total degree of polynomial}.

Hereafter, for cases $n=2$ and $n=3$, we denote $(k_1,k_2)=:(k,l)$ and $(k_1,k_2,k_3)=:(k,l,m)$.

The level vectors for Lie algebras of ranks two and three are in the Appendix, for all cases see~\cite{McP}.

\begin{example}
Consider group $B_3$. Its level vector equals $(6,10,6)$.
Let us order monomials $u_1^2$, $u_2^2$ and $u_3^3$ in the case of $B_3$.
To do this we respectively calculate vectors $\mathfrak{n}$, $\mathfrak{n}'$ and $\mathfrak{n}''$:
\begin{gather*}
\mathfrak{n}=(6,10,6)(2,0,0)^t,\quad
\mathfrak{n}'= (6,10,6)(0,2,0)^t,\quad
\mathfrak{n}''= (6,10,6)(0,0,3)^t.
\end{gather*}
Therefore we obtain $u_2^2\succ u_3^3\succ u_1^2$.

Similarly for the few first degrees we have:
\begin{gather*}
u_2^2 \succ u_1^3\succ u_1^2u_3\succ u_1u_3^2\succ u_3^3\succ u_1u_2\succ u_2u_3\succ u_1^2\succ u_1u_3\succ u_3^2\succ u_2\succ u_1\succ u_3\succ 1.
\end{gather*}
\end{example}

\begin{definition}
Let $\{P_{k_1,\ldots, k_n}(u)\} \in \C[u_1,\ldots,u_n]$ be a family of polynomials in $n$ variables satisfying
\begin{gather*}
\int\limits_{\C^n} P_{k_1,\ldots ,k_n}(x)P_{k_1',\ldots, k_n'}(x)d \rho(x)=
\sum_{u\in V}P_{k_1,\ldots, k_n}(u)P_{k_1',\ldots, k_n'}(u)\varrho(u)=
\Lambda_{k_1,\ldots, k_n}\delta_{k_1k_1'}\cdots \delta_{k_nk_n'},
\end{gather*}
where $d \rho (x) = \sum\limits_{u\in V}\varrho(x)\delta(x-u)$ is the discrete measure,  $V$ is a lattice in $\C^n$
and $\Lambda_{k_1,\ldots, k_n}$ is the normalization constant (see~\cite{Dunkl}, \cite{Xu} and \cite{IlievXu}).

The family $\{P_{k}(u)\}$ is then called \textit{orthogonal polynomials} in $n$ variables.
\end{definition}

The orthogonal polynomials have a number of useful properties, in particular each polynomial satisfies recurrence relation,
see, e.g.~\cite{Xu} and~\cite{IlievXu}.

\subsection{Orthogonal polynomials of simple group $G$}
Here we use the approach to the construction of orthogonal polynomials in $n$ variables that was proposed in~\cite{NPT}.
It is based on the idea of replacement of the lowest $C$-orbit functions $C_{\w_j}$, $j=1,2,\dots,n$ by new variables $X_j$.
This method brought us to the generalization of classical Chebyshev polynomials
to Chebyshev polynomials in $n$-dimensional Euclidean~space.
Moreover it easily gives us rather wide families of Macdonald polynomials.
Polynomials generated from $C$-functions can be viewed as the generalized Chebyshev polynomials of the first kind in $n$ variables
and as the Macdonald symmetric polynomials for the case $k_\alpha=0$, $t_\alpha=1$.
$S$-polynomials play role of generalized Chebyshev polynomials of the second kind
and equivalent to the Macdonald polynomials with $k_\alpha=1$, $t_\alpha=q_\alpha$.
The $S$-functions divided by the lowest $S$-function $S_\rho(x)$ coincide with the character of the representation
and we do use these fractions as the $S$-polynomials.


Consider $C$-functions and $S$-functions defined in the Preliminaries and introduce new coordinates 
\begin{gather}\label{var_subst}
u_1:=C_{(1,0,\ldots, 0)}(x),\quad u_2:=C_{(0,1,0,\ldots, 0)}(x),\; \dots,\; u_n:=C_{(0,\ldots, 0,1)}(x).
\end{gather}

These variables coincide with those introduced in~\cite{NPT}: $X_j=C_{\w_j}=u_j$, $j=1,2,\dots,n$.
It was shown in~\cite{NPT} (see Proposition 1) that in these coordinates orbit functions gain the polynomial form,
and it directly follows from the orthogonality of orbit functions that polynomials $C_{\lambda}(u)$ and $S_{\lambda}(u)$ are orthogonal.
The orthogonality regions for polynomials $C_{\lambda}(u)$ and $S_{\lambda}(u)$ are the images $\mathfrak{F}$ and $\widetilde{\mathfrak{F}}$
of the fundamental regions $F$ and $\widetilde{F}$ under the transformation $x=(x_1,\dots,x_n)\mapsto u=(u_1,\dots,u_n)$.
The discretization developed for the orbit functions (see Section~\ref{sec_discretization})
can be effectively applied to these polynomials.
Let us fix positive $M$ and let $x=(\frac{s_1}{M},\ldots,\frac{s_n}{M})\in F_M$. Using (\ref{Corth}) we have
\begin{gather*}
\langle C_\lambda(x), C_{\lambda'}(x) \rangle_{F_M}{=}
\!\!\sum\limits_{x\in {F_M}}\!{|W_x|} C_\lambda(x)\overline{C_{\lambda'}(x)}{=}
\!\!\sum\limits_{u\in {\mathfrak{F}_M}}\!|W_u|J(u) C_\lambda(u)\overline{C_{\lambda'}(u)}{=}
\det \mathfrak{C}\ \tfrac{|W|^2M^n}{|W_\lambda|}\delta_{\lambda\ \lambda'},
\end{gather*}
where $J^{-1}(u)$ 
is the discretized Jacobian of the transformation $x\mapsto u$.

Thereby we obtain discrete orthogonality of polynomials $C_{\lambda}(u)$ with the weight function~${|W_u|}J(u)$.

Let us calculate $J(u)$ explicitly.
Whereas $\rho{=}\w_1{+}\dots{+}\w_n=(1,\dots,1)_{\w}$ and $\mathbf{S}(u):=S_{\rho}^2(x)$, it is easy to check that
Jacobian $J(u)=\frac{1}{(2\pi)^n  \sqrt{|{\mathbf{S}(u)}|}}$\quad for $A_n$, $B_n$ and $C_n$.
Similarly for the $S$-polynomials defined as character $\mathcal{S}_{\lambda}(u):=\frac{S_{\lambda+\rho}(x)}{S_{\rho}(x)}$ from~(\ref{Sorth}), we have
the discrete orthogonality of characters
\begin{gather*}
\langle S_\lambda(x), S_{\lambda'}(x) \rangle_{\widetilde{F}_M}{=}
\!\!\!\sum\limits_{x\in\widetilde{F}_M}\!\!{|W|}S_\lambda(x)\overline{S_{\lambda'}(x)}{=}
\!\!\!\sum\limits_{u\in{\widetilde{\mathfrak{F}}_M}}\!\!|W|J(u){\mathbf{S}(u)}\mathcal{S}_\lambda(u)\overline{\mathcal{S}_{\lambda'}(u)}{=}
\det \mathfrak{C}|W| M^n\delta_{\lambda\,\lambda'}.
\end{gather*}

Hereof we obtained discrete orthogonality of polynomials  with the weight function~${|W|}J(u){\mathbf{S}(u)}$ and\quad
$J(u)=\frac{1}{(2\pi)^n\sqrt{|{\mathbf{S}(u)}|}}$\quad for $A_n$, $B_n$ and $C_n$.

Explicit values of the weight functions for the groups
$A_2$, $C_2$, $A_3$, $B_3$ and $C_3$ are presented in Appendix.

\section{Polynomials of two complex variables}
This section recalls recursion relations and lowest orthogonal polynomials generated from groups $A_2$ and $C_2$.
We need them here in order to recognize the result of the reduction to the maximal subgroup in Section~\ref{sec_reductions}.
In contrast to the paper~\cite{NPT}, we consider explicit examples of discretization of polynomials
and transformed fundamental regions $\mathfrak{F}$.

In this section, we mean that $S$-polynomials are the fractions $\frac{S_{(k+1,l+1)}(x)}{S_{(1,1)}(x)}$,
and we denote them by  $\mathcal{S}_{(k,l)}(u)$.

\subsection{$C$- and $S$-polynomials of $A_2$}
Let us introduce new coordinates $u :=(u_1,u_2)=(C_{(1,0)}(x), C_{(0,1)}(x))$.
Using the decomposition of products of orbit functions (see~\cite{NPT}) we obtain the set of recurrence relations.
In the case of $S$-polynomials we adduce only generic recurrence relations, and the lowest $S$-polynomials
can be obtained from $C$-polynomials using the Weyl character formula and multiplicities from~\cite{BMP}.
The lowest $C$-polynomials are constructed and arranged in Table~\ref{c_poly_A2}.

Generic recurrence relations for $S$-polynomials of $A_2$ (Chebyshev polynomials of the second kind in two variables):
\begin{gather*}
\mathcal{S}_{(k + 1,l)}(u)=u_1 \mathcal{S}_{(k,l)}(u) - \mathcal{S}_{(k,l-1)}(u) - \mathcal{S}_{(k-1,l+1)}(u),\qquad k,l>1;\\
\mathcal{S}_{(k,l + 1)}(u)=u_2 \mathcal{S}_{(k,l)}(u) - \mathcal{S}_{(k-1,l)}(u)-\mathcal{S}_{(k + 1,l-1)}(u),\qquad k,l>1.
\end{gather*}

Recurrence relations for $C$-polynomials of $A_2$ (Chebyshev polynomials of the first kind in two variables):
\begin{gather*}
C_{(k+1,l)}(u)=u_1C_{(k,l)}(u)-C_{(k,l-1)}(u)-C_{(k-1,l+1)}(u),\qquad k,l>1;\\
C_{(k,l+1)}(u)=u_2C_{(k,l)}(u)-C_{(k+1,l-1)}(u)-C_{(k-1,l)}(u),\qquad k,l>1;\\
C_{(k+1,0)}(u)=u_1C_{(k,0)}(u)-C_{(k-1,1)}(u),\qquad k>1;\\
C_{(0,l+1)}(u)=u_2C_{(0,l)}(u)-C_{(1,l-1)}(u),\qquad l>1.
\end{gather*}
\begin{table}[h]
\begin{tabular}{ccc}
$\sharp 0$  &  $\sharp 1$   &   $\sharp 2$\\
$\begin{array}{|l|r|r|r|r|}
\hline\hline
C_{(k, l)}(u)&1 &u_1u_2&u_1^3&u_2^3\\ \hline
C_{(0,0)}(u) &1 &      &     &     \\ \hline
C_{(1,1)}(u) &-3&1     &     &     \\ \hline
C_{(3,0)}(u) &3 &-3    & 1   &     \\ \hline
C_{(0,3)}(u) &3 &-3    & 0   &1    \\ \hline
\end{array}$
&
$\begin{array}{|l|r|r|r|}
\hline \hline C_{(k, l)}(u)&u_2& u_1^2& u_1u_2^2   \\ \hline
C_{(0,1)}(u) & 1 &       &               \\ \hline
C_{(2,0)}(u) &-2 & 1     &             \\ \hline
   C_{(1,2)}(u) & -1 &-2     &1 \\ \hline
       \end{array}    $    &
$ \begin{array}{|l|r|r|r|}
   \hline \hline C_{(k, l)}(u)&u_1& u_2^2&u_1^2u_2 \\ \hline
    C_{(1,0)}(u)  & 1 &      &       \\ \hline
    C_{(0,2)}(u)  &-2 & 1    &       \\ \hline
    C_{(2,1)}(u)  &-1 & -2    & 1     \\ \hline
                    \end{array}$
\end{tabular}
\caption{Lowest $C_{(k,l)}$-polynomials of the Lie group $A_2$ split into three classes
in correspondence with the congruence number $\#$ of the parameters $(k,l)$.}\label{c_poly_A2}
\end{table}

Let us study the transformation~\eqref{var_subst} of the fundamental region $F(A_2)\to \mathfrak{F}(A_2)$.
During the substitution $x\mapsto u$, the vertices of the simplex $F(A_2)$ (see Appendix) go to the points $\{P_0,P_1,P_2\}$ and edges  go to the continuous curves.
Explicitly we have
\begin{gather*}
(0,0) \mapsto (3,0)=:P_0;\qquad
\omega_1 \mapsto (-\tfrac{3}{2},-\tfrac{3\sqrt 3}{2})=:P_1;\qquad
\omega_2 \mapsto (-\tfrac{3}{2},\tfrac{3\sqrt 3}{2})=:P_2.
\end{gather*}

\begin{example}\label{ex_a2}
Let us fix $M=3$. $|F_3(A_2)|=10$ and the corresponding grid points $(\frac{s_1}{M},\frac{s_2}{M})$ in coordinates
$(\textrm{Re}(u_1),\textrm{Im}(u_1))$ are
\begin{gather*}
\begin{array}{@{\!\!\!\!\!\!\!}ll}
(0,0)                      \mapsto(3,0);                                                                                                       &\;  (0,1)\mapsto({-}\tfrac{3}{2},\tfrac{3\sqrt 3}{2});\\
(\tfrac{2}{3},\tfrac{1}{3})\mapsto({-}2 \cos{\tfrac{\pi}{9}}{+} \cos{\tfrac{2 \pi}{9}}, {-}2 \sin{\tfrac{\pi}{9}} {-} \sin{\tfrac{2 \pi}{9}}); &\;  (1,0)\mapsto({-}\tfrac{3}{2},{-}\tfrac{3\sqrt 3}{2});\\
(\tfrac{2}{3},0)           \mapsto({-} \cos{\tfrac{\pi}{9}} {+} 2\sin{\tfrac{ \pi}{18}}, {-}2 \cos{\tfrac{\pi}{18}} {+} \sin{\tfrac{ \pi}{9}});&\; (\tfrac{1}{3},\tfrac{2}{3})\mapsto({-}2 \cos{\tfrac{\pi}{9}} {+} \cos{\tfrac{2 \pi}{9}}, 2 \sin{\tfrac{\pi}{9}} {+} \sin{\tfrac{2 \pi}{9}});\\
(\tfrac{1}{3},0)           \mapsto(2 \cos{\tfrac{2\pi}{9}}  {+} \sin{\tfrac{ \pi}{18}},  \cos{\tfrac{\pi}{18}} {-} 2\sin{\tfrac{2 \pi}{9}});   &\; (0,\tfrac{1}{3})           \mapsto(2 \cos{\tfrac{2\pi}{9}}   {+} \sin{\tfrac{ \pi}{18}}, {-} \cos{\tfrac{\pi}{18}} {+}2 \sin{\tfrac{2\pi}{9}});\\
(\tfrac{1}{3},\tfrac{1}{3})\mapsto(0,0);                                                                                                       &\; (0,\tfrac{2}{3})\mapsto({-}2 \cos{\tfrac{\pi}{9}} {+} 2\sin{\tfrac{ \pi}{18}}, 2 \cos{\tfrac{\pi}{18}} {-} \sin{\tfrac{ \pi}{9}}).
\end{array}
\end{gather*}
The choice of new coordinates in the form $(\textrm{Re}(u_1),\textrm{Im}(u_1))$ is determined by the complex conjugation $u_1=\overline{u_2}$.
Plotting these points and the transformed fundamental region $\mathfrak{F}$ in Figure~\ref{fig_a2},
we see that the discrete set $F_M(A_2)$ in new coordinates $(u_1,u_2)$ goes to points $\mathfrak{F}_M(A_2)\subset\mathfrak{F}(A_2).$
\begin{figure}[h]
\centerline{\includegraphics[scale=0.7]{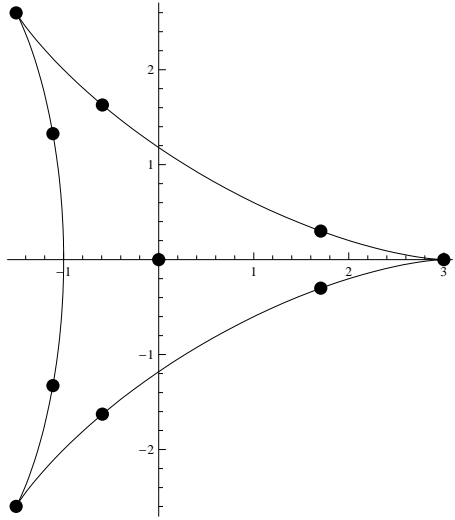}}
\caption{Region of orthogonality $\mathfrak{F}$ of polynomials of $A_2$ and discrete points of $\mathfrak{F}_3$ obtained in Example~\ref{ex_a2}.}\label{fig_a2}
\end{figure}
\end{example}

\subsection{$C$- and $S$-polynomials of $C_2$}
As in the previous section, we obtain by the same substitution the orthogonal polynomials of group~$C_2$.
Recursion relations for these $S$- and $C$-polynomials are also obtained by the expansion of the products of orbit functions.

The generic recurrence relations for $C$-polynomials are the following
\begin{gather*}
C_{(k+1,l)}(u)=u_1C_{(k,l)}(u)-C_{(k-1,l)}(u)-C_{(k+1,l-1)}(u)-C_{(k-1,l+1)}(u),\quad k,l>1;\\
C_{(k,l+1)}(u)=u_2C_{(k,l)}(u)-C_{(k,l-1)}(u)-C_{(k+2,l-1)}(u)-C_{(k-2,l+1)}(u),\quad k>2,\ l>1.\\
\end{gather*}
Some of the additional relations are:
\begin{gather*}
C_{(k+1,0)}(u)=u_1C_{(k,0)}(u)-C_{(k-1,0)}(u)-C_{(k-1,1)}(u),\quad k>1;\\
C_{(0,l+1)}(u)=u_2C_{(0,l)}(u)-C_{(0,l-1)}(u)-C_{(2,l-1)}(u),\quad l>1.
\end{gather*}
The remaining low-order $C$-polynomials necessary to solve all above recursions are presented in Table~\ref{c_poly_A2}.

When we have the lowest $C$-polynomials, all $S$-polynomials can be found from the relations:
\begin{gather*}
\mathcal{S}_{(k+1,l)}(u)=u_1\mathcal{S}_{(k,l)}(u)-\mathcal{S}_{(k+1,l-1)}(u)-\mathcal{S}_{(k-1,l+1}(u)-\mathcal{S}_{(k-1,l)}(u),\quad k,l>1;\\
\mathcal{S}_{k\l+1}(u)=u_2\mathcal{S}_{(k,l)}(u)-\mathcal{S}_{(k,l-1)}(u)-\mathcal{S}_{(k+2,l-1)}(u)-\mathcal{S}_{(k-2,l+1)}(u),\quad k>2,\ l>1.
\end{gather*}
\begin{table}[h]
\begin{tabular}{cc}
$\begin{array}{|l|r|r|r|r|r|r|}
\hline\hline
C_{(k,l)}(u)&1&u_2&u_1^2&u_2^2&u_1^2u_2&u_2^3\\\hline
C_{(0,0)}(u)&1&&&&&\\\hline
C_{(0,1)}(u)&0&1&&&&\\\hline
C_{(2,0)}(u)&-4&-2&1&&&\\\hline
C_{(0,2)}(u)&4&4&-2&1&&\\\hline
C_{(2,1)}(u)&0&-6&0&-2&1&\\\hline
C_{(0,3)}(u)&0&9&0&6&-3&1\\\hline
\end{array}$
&\qquad
$\begin{array}{|l|r|r|r|r|}
\hline\hline
C_{(k,l)}(u)&u_1&u_1u_2&u_1^3&u_1u_2^2\\\hline
C_{(1,0)}(u)&1&&&\\\hline
C_{(1,1)}(u)&-2&1&&\\\hline
C_{(3,0)}(u)&-3&-3&1&\\\hline
C_{(1,2)}(u)&6&3&-2&1\\\hline
\end{array}$
\end{tabular}
\caption{Lowest $C$-polynomials of $C_2$ split into two congruence classes
according to the parity of the parameter $k$.}\label{c_poly_c2}
\end{table}

Under the transformation $x\mapsto u$, the vertices of the simplex $F(C_2)$ (see Appendix) go to points $\{P_0,P_1,P_2\}$,
namely\quad
$
(0,0)    \mapsto (4,4 )=:P_0,\quad
\omega_1 \mapsto (0,-4)=:P_1,\quad
\omega_2 \mapsto (-4,4)=:P_2.
$

\begin{example}
Let us fix $M=4$. $|F_4(C_2)|=9$ and the corresponding grid points $(\frac{s_1}{M},\frac{s_2}{M})$ in coordinates
$(u_1,u_2)$ are:
\begin{gather*}
\begin{array}{lll}
(0,0)           \mapsto (4,4),       & (0,\frac{1}{4}) \mapsto (2\sqrt2,2),       & (\frac{1}{2},0)          \mapsto (0,-4),\\
(0,\frac{1}{2}) \mapsto (0,0),       & (0,\frac{3}{4}) \mapsto (-2\sqrt2,2),\quad & (\frac{1}{4},\frac{1}{2})\mapsto (-2,0),\\
(0,1)           \mapsto (-4,4),\quad & (\frac{1}{4},0) \mapsto (2,0),             &(\frac{1}{4},\frac{1}{4}) \mapsto (0,-2).
\end{array}
\end{gather*}

In Figure~\ref{fig_c2} we see the discrete set $F_M(C_2)$ in new coordinates $(u_1,u_2)$ and new fundamental region $\mathfrak{F}(C_2).$
\begin{figure}[h]
\centerline{\includegraphics[scale=0.8]{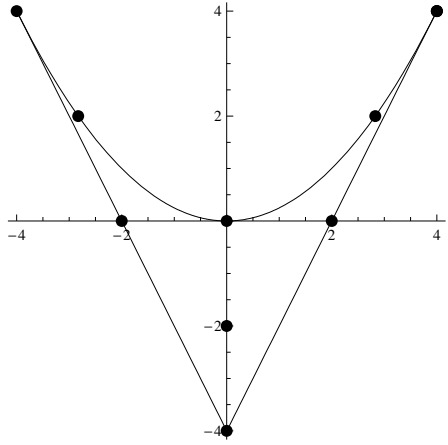}}
\caption{Region of orthogonality $\mathfrak{F}$ and discrete points $\mathfrak{F}_4$ of polynomials of $C_2$.}\label{fig_c2}
\end{figure}
\end{example}

\section{Polynomials in three variables}
In this section, we first recall recursion relations and orthogonal polynomials of group $A_3$ and
then we present new orthogonal polynomials of groups $B_3$ and $C_3$
together with the complete sets of recurrence relations.
The orthogonality domains of the polynomials of all these groups and discretization examples are shown.

For the orbit functions of each of the groups $A_3$, $B_3$ and $C_3$, we introduce the new coordinates by the rule
$x\mapsto u$, where $x=(x_1,x_2,x_3)$, $u=(u_1, u_2, u_3)$ and
\begin{gather}\label{3d_subst}
u_1:=C_{\w_1}(x)=C_{(1,0,0)}(x),\quad
u_2:=C_{\w_2}(x)=C_{(0,1,0)}(x),\quad
u_3:=C_{\w_3}(x)=C_{(0,0,1)}(x).
\end{gather}

Note that for different groups of rank three, the lowest $C$-functions $u_1$, $u_2$, $u_3$ have different forms and properties,
but the general substitution rule~\eqref{3d_subst} is always the same.

The recurrence relations for $C$- and $S$- polynomials come from the expansions of the products
$u_jC_{\lambda}$ and $u_j\mathcal{S}_{\lambda}$, $j=1,2,3$.

We also mean that $S$-polynomials are the fractions $\frac{S_{\lambda+\rho}(x)}{S_{\rho}(x)}$
and we denote them as  $\mathcal{S}_{\lambda}(u)$.

As in this section $\lambda$ has only three coordinates in $\w$-basis, we use the notations
$\lambda=(k,l,m)$, $C_{(k,l,m)}(u)$ and $\mathcal{S}_{(k,l,m)}(u)$.

\subsection{$C$- and $S$-polynomials of $A_3$ in three complex variables}
When the lowest $C$-functions $C_{\w_1}(x)$, $C_{\w_2}(x)$ and $C_{\w_3}(x)$ of $A_3$ are written explicitly,
it is easy to verify that $u_2=\overline{u}_2$ and $u_3=\overline{u}_1$.

Generic recursions for $C$-polynomials (Chebyshev polynomials of the first kind) in three complex variables are:
\begin{gather*}
C_{(k{+}1,l,m)}(u)=u_1C_{(k,l,m)}(u){-}C_{(k{-}1,l{+}1,m)}(u){-}C_{(k,l,m{-}1)}(u){-}C_{(k,l{-}1,m{+}1)}(u),\;\; k,l,m\geqslant 1;\\
C_{(k,l{+}1,m)}(u)=u_2C_{(k,l,m)}(u){-}C_{(k{+}1,l,m{-}1)}(u){-}C_{(k{+}1,l{-}1,m{+}1)}(u)\\
\phantom{C_{(k,l{+}1,m)}(u)=}       {-}C_{(k{-}1,l{+}1,m{-}1)}(u){-}C_{(k{-}1,l,m{+}1)}(u){-}C_{(k,l{-}1,m)}(u),\;\; k,l,m\geqslant 1;\\
C_{(k,l,m{+}1)}(u)=u_3C_{(k,l,m)}(u){-}C_{(k{+}1,l{-}1,m)}(u){-}C_{(k{-}1,l,m)}(u){-}C_{(k,l{+}1,m{-}1)}(u),\;\; k,l,m\geqslant 1.
\end{gather*}

Note that the last of the above relations can also be obtained from
$C_{(k,l,m)}(u) = \overline{C_{(m,l,k)}(u)}$ for $k,l,m \in \Z^{\geqslant 0}.$

Additional recursion relations for  $C$-polynomials of $A_3$:
\begin{gather*}
C_{(k+1,0,0)}(u)=u_1C_{(k,0,0)}(u)-C_{(k-1,1,0)}(u),\quad k\geqslant 1;\\
C_{(0,k+1,0)}(u)=u_2C_{(0,k,0)}(u)-C_{(0,k-1,0)}(u)-C_{(1,k-1,1)}(u),\quad k\geqslant 1;\\
C_{(k+1,l,0)}(u)=u_1C_{(k,l,0)}(u)-C_{(k,l- 1,1)}(u)-C_{(k-1,l+ 1,0)}(u),\quad k,l\geqslant 1;\\
C_{(k,l+1,0)}(u)=u_2C_{(k,l,0)}(u)-C_{(k-1,l,1)}(u)-C_{(k,l-1,0)}(u)-C_{(k+1,l-1,1)}(u),\quad k,l\geqslant 1;\\
C_{(k+1,0,m)}(u)=u_1C_{(k,0,m)}(u)-C_{(k-1,1,m)}(u)-C_{(k,0,m-1)}(u)\quad k,m\geqslant 1.
\end{gather*}

The remaining low-order $C$-polynomials of $A_3$ can be found in Table~\ref{c_poly_a3}.

\begin{table}[h]
$\begin{array}{cc}
\sharp 0                               & \sharp 1\\
\begin{array}{|l|r|r|r|r|r|}
\hline\hline
C_{(k, l,m)}(u)&1&  u_1u_3 & u_2^2&u_2u_3^2&u_1^2u_2 \\ \hline
C_{(0,0,0)}(u) &1&        &       &       &        \\ \hline
C_{(1,0,1)}(u) &-4& 1     &       &       &        \\ \hline
C_{(0,2,0)}(u) &2&  -2    & 1     &       &        \\ \hline
C_{(0,1,2)}(u) &4&  -1    & -2    &   1   &        \\ \hline
C_{(2,1,0)}(u) &4&  -1    & -2    &   0   &    1    \\ \hline
\end{array}
\quad&
\begin{array}{|l|r|r|r|r|r|}
\hline \hline
C_{(k, l,m)}(u)&u_1& u_2u_3& u_3^3&u_1^2u_3&u_1u_2^2   \\ \hline
C_{(1,0,0)}(u) & 1 &       &      &       &         \\ \hline
C_{(0,1,1)}(u) &-3 & 1     &      &       &        \\ \hline
C_{(0,0,3)}(u) & 3 &-3     &1     &       &        \\ \hline
C_{(2,0,1)}(u) &-1 &-2     & 0    &   1   &        \\ \hline
C_{(1,2,0)}(u) & 5 & -1    & 0    &  -2   &   1    \\ \hline
\end{array}
\end{array}$
\bigskip

$\begin{array}{cc}
\sharp 2 & \sharp 3
\\
\begin{array}{|l|r|r|r|r|r|}
\hline \hline
C_{(k, l,m)}(u) &u_2& u_3^2&u_1^2 &u_1u_2u_3&u_2^3  \\ \hline
C_{(0,1,0)}(u)  & 1 &      &      &       & \\ \hline
C_{(0,0,2)}(u)  &-2 & 1    &      &       & \\ \hline
C_{(2,0,0)}(u)  &-2 & 0    & 1    &       & \\ \hline
C_{(1,1,1)}(u)  & 4 & -3   &-3    &  1    & \\ \hline
C_{(0,3,0)}(u)  & -3&  3   & 3    & -3    &1\\ \hline
\end{array}
\qquad&
\begin{array}{|l|r|r|r|r|r|}
\hline \hline
C_{(k, l,m)}(u)& u_3& u_1u_2& u_1^3&u_1u_3^2&u_2^2u_3 \\ \hline
C_{(0,0,1)}(u) &  1 &       &      &        &      \\ \hline
C_{(1,1,0)}(u) & -3 &   1   &      &        &      \\ \hline
C_{(3,0,0)}(u) &  3 &  -3   &  1   &        &       \\ \hline
C_{(1,0,2)}(u) & -1 & -2    & 0    &  1     &      \\ \hline
C_{(0,2,1)}(u) &  5 & -1    & 0    &  -2    & 1    \\ \hline
\end{array}
\end{array}$
\caption{Lowest $C$-polynomials of $A_3$ split into four congruence classes
$\#=0$, $\#=1$, $\#=2$ and $\#=3$.}\label{c_poly_a3}
\end{table}

Let us study the transformation of the fundamental region $F(A_3)\to \mathfrak{F}(A_3)$.
During the substitution $x\mapsto u$ the vertices of the simplex $F(A_3)$ (see Appendix) go to the points $\{P_0,P_1,P_2,P_3\}$:
\begin{gather*}
(0,0,0)  \mapsto  (4,6,0) =:P_0,\qquad
\omega_1 \mapsto  (0,{-}6,{-}4)=:P_1,\\
\omega_2 \mapsto  ({-}4,6,0 )=:P_2,\qquad\qquad
\omega_3 \mapsto  (0,{-}6,4)=:P_3.
\end{gather*}
The shape of the region of orthogonality of polynomials of $A_3$ is presented in Figure~\ref{fig_a3}.


\begin{example}
Let us consider discretization with $M=3$. $|F_3(A_3)|=20$ and grid points $(\frac{s_1}{M},\frac{s_2}{M},\frac{s_3}{M})$
in new coordinates $(\textrm{Re}(u_1),u_2,\textrm{Im}(u_1))$ are:

\begin{gather*}
\begin{array}{lll}
(0,0,0)                                 \mapsto (4,6,0),                                     & (0,0,1)                       \mapsto (0,-6,4),                                     & (0,1,0)                       \mapsto (-4,6,0),                             \\
(\tfrac{2}{3},0,\tfrac{1}{3})           \mapsto (0,-3,-2),                                   & (\tfrac{2}{3},\tfrac{1}{3},0) \mapsto (-\tfrac{1}{2},-3,-\tfrac{3\sqrt3}{2}),\qquad & (\tfrac{1}{3},0,\tfrac{2}{3}) \mapsto (0,-3,2),                             \\
(\tfrac{1}{3},\tfrac{2}{3},0)           \mapsto (-\tfrac{3\sqrt3}{2},3,-\tfrac{1}{2}),\qquad & (0,\tfrac{2}{3},\tfrac{1}{3}) \mapsto (-\tfrac{3\sqrt3}{2},3,\tfrac{1}{2}),         & (\tfrac{1}{3},0,0)            \mapsto (\tfrac{3\sqrt3}{2},3,-\tfrac{1}{2}), \\
(0,\tfrac{1}{3},0)                      \mapsto (2,3,0),                                     & (\tfrac{1}{3},\tfrac{1}{3},0) \mapsto (0,0,-1),                                     & (\tfrac{1}{3},0,\tfrac{1}{3}) \mapsto (1,0,0),                              \\
(\tfrac{1}{3},\tfrac{1}{3},\tfrac{1}{3})\mapsto (-1,0,0),                                    & (\tfrac{2}{3},0,0)            \mapsto (\tfrac{1}{2},-3,-\tfrac{3\sqrt3}{2}),        & (0,0,\tfrac{2}{3})            \mapsto (\tfrac{1}{2},-3,\tfrac{3\sqrt3}{2}), \\
(1,0,0)                                 \mapsto (0,-6,-4),                                   & (0,\tfrac{1}{3},\tfrac{2}{3}) \mapsto (-\tfrac{1}{2},-3,\tfrac{3\sqrt3}{2}),        & (0,0,\tfrac{1}{3})            \mapsto (\tfrac{3\sqrt3}{2},3,\tfrac{1}{2}),  \\
(0,\tfrac{1}{3},\tfrac{1}{3})           \mapsto (0,0,1),                                     & (0,\tfrac{2}{3},0)            \mapsto (-2,3,0).&
\end{array}
\end{gather*}
All these points belong to the new fundamental region $\mathfrak{F}(A_3)$.

\begin{figure}[h]
\centerline{\includegraphics[scale=0.8]{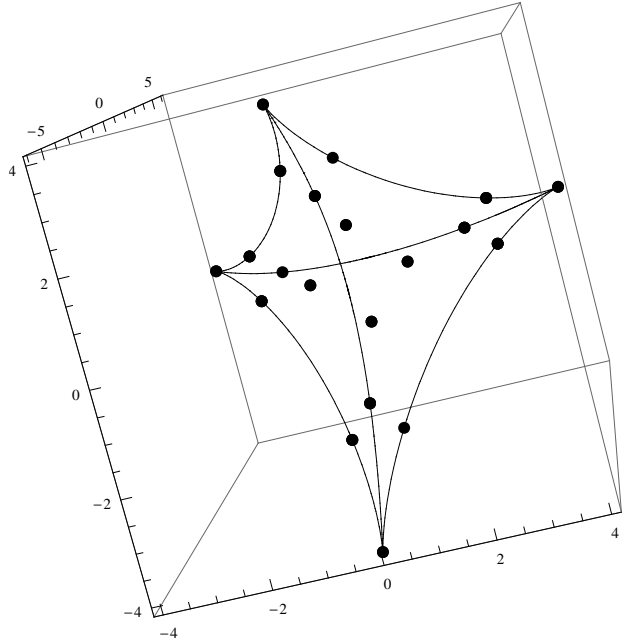}}
\caption{Region of orthogonality $\mathfrak{F}$ and discrete points $\mathfrak{F}_3$ of polynomials of $A_3$.}\label{fig_a3}
\end{figure}
\end{example}

Using the Weyl character formula and multiplicities from~\cite{BMP}, we can obtain $S$-polynomials
(Chebyshev polynomials of the second kind) in $u_1$, $u_2$ and $u_3$.
Or, alternately, having the lowest $S$-polynomials we can construct other polynomials
by means of the following recursions
\begin{gather*}
\mathcal{S}_{(k{+}1,l,m)}(u)=u_1\mathcal{S}_{(k,l,m)}(u){-}\mathcal{S}_{(k{-}1,l{+}1,m)}(u){-}\mathcal{S}_{(k,l,m{-}1)}(u){-}\mathcal{S}_{(k,l{-}1,m{+}1)}(u),\quad k,l,m\geqslant 2;\\
\mathcal{S}_{(k,l + 1,m)}(u)=u_2\mathcal{S}_{(k,l,m)}(u){-}\mathcal{S}_{(k,l{-}1,m)}(u){-}\mathcal{S}_{(k{+}1,l{-}1,m{+}1)}(u){-}\mathcal{S}_{(k{-}1,l{+}1,m{-}1)}(u)\\
\phantom{\mathcal{S}_{(k,l{+}1,m)}(u)=}     {-}\mathcal{S}_{(k{+}1,l,m{-}1)}(u){-}\mathcal{S}_{(k{-}1,l,m{+}1)}(u),\quad k,l,m\geqslant 2;\\
\mathcal{S}_{(k,l,m{+}1)}(u)=u_3\mathcal{S}_{(k,l,m)}(u){-}\mathcal{S}_{(k{+}1,l{-}1,m)}(u){-}\mathcal{S}_{(k{-}1,l,m)}(u){-}\mathcal{S}_{(k,l{+}1,m{-}1)}(u),\quad k,l,m\geqslant 2.
\end{gather*}

\subsection{$C$- and $S$-polynomials of $B_3$ in three real variables}
For group $B_3$, our new coordinates $u$ satisfy the relation $u_i=\overline{u}_i$, $i=1,2,3$.

The following generic recursion relations for $C$-polynomials hold true when $k,l,m\geqslant 2$:
\begin{gather*}
C_{(k{+}1,l,m)}(u)=u_1C_{(k,l,m)}(u){-}C_{(k,l{+}1,m{-}2)}(u){-}C_{(k,l{-}1,m{+}2)}(u){-}C_{(k{+}1,l{-}1,m)}(u)\\
\phantom{C_{(k{+}1,l,m)}(u)=}                    {-}C_{(k{-}1,l{+}1,m)}(u){-}C_{(k{-}1,l,m)}(u){-}C_{(k{-}1,l{+}1,m)}(u){-}C_{(k{-}1,l,m)}(u),\\
C_{(k,l{+}1,m)}(u)=u_2C_{(k,l,m)}(u){-}C_{(k{+}1,l{-}1,m{+}2)}(u){-}C_{(k{-}1,l{-}1,m{+}2)}(u){-}C_{(k{+}1,l{-}2,m{+}2)}(u)\\
\phantom{C_{(k,l{+}1,m)}(u)=}                    {-}C_{(k{-}1,l,m{+}2)}(u){-}C_{(k{-}1,l{+}2,m{-}2)}(u){-}C_{(k{+}1,l{+}1,m{-}2)}(u)\\
\phantom{C_{(k,l{+}1,m)}(u)=}                    {-}C_{(k{-}1,l{+}1,m{-}2)}(u){-}C_{(k{+}1,l,m{-}2)}(u){-}C_{(k{-}2,l{+}1,m)}(u){-}C_{(k{+}2,l{-}1,m)}(u),\\
C_{(k,l,m{+}1)}(u)=u_3C_{(k,l,m)}(u){-}C_{(k{+}1,l{-}1,m{+}1)}(u){-}C_{(k,l{-}1,m{+}1)}(u){-}C_{(k{-}1,l,m{+}1)}(u)\\
\phantom{C_{(k,l,m{+}1)}(u)=}                    {-}C_{(k{-}1,l{+}1,m{-}1)}(u){-}C_{(k,l{+}1,m{-}1)}(u){-}C_{(k{+}1,l,m{-}1)}(u){-}C_{(k,l,m{-}1)}(u)
\end{gather*}

Remaining recurrence relations except for the lowest polynomials are listed below:
\begin{gather*}
C_{(k{+}1,l,0)}(u)=u_1C_{(k,l,0)}(u){-}C_{(k{-}1,l{+}1,)0}(u){-}C_{(k{+}1,l{-}1,0)}(u){-}C_{(k{-}1,l,0)}(u){-}C_{(k,l{-}1,2)}(u),\\
C_{(k,l{+}1,0)}(u)=u_2C_{(k,l,0)}(u){-}C_{(k{-}2,l{+}1,0)}(u){-}C_{(k{+}2,l{-}1,0)}(u){-}C_{(k{-}1,l,2)}(u){-}C_{(k,l{-}1,0)}(u)\\
\phantom{C_{(k,l{+}1,0)}(u)=}                    {-}C_{(k{+}1,l{-}1,2)}(u){-}C_{(k{-}1,l{-}1,2)}(u){-}C_{(k{+}1,l{-}2,2)}(u),\\
C_{(k{+}1,0,0)}(u)=u_1C_{(k,0,0)}(u){-}C_{(k{-}1,1,0)}(u){-}C_{(k{-}1,0,0)}(u),\\
C_{(0,l{+}1,0)}(u)=u_2C_{(0,l,0)}(u){-}C_{(2,l{-}1,0)}(u){-}C_{(0,l{-}1,0)}(u){-}C_{(1,l{-}1,2)}(u){-}C_{(1,l{-}2,2)}(u),\\
C_{(0,0,m{+}1)}(u)=u_3C_{(0,0,m)}(u){-}C_{(0,1,m{-}1)}(u){-}C_{(1,0,m{-}1)}(u){-}C_{(0,0,m{-}1)}(u),\\
C_{(k{+}1,0,m)}(u)=u_1C_{(k,0,m)}(u){-}C_{(k,1,m{-}2)}(u){-}C_{(k{-}1,1,m)}(u),\\
C_{(k,0,m{+}1)}(u)=u_3C_{(k,0,m)}(u){-}C_{(k{-}1,0,m{+}1)}(u){-}C_{(k{-}1,1,m{-}1)}(u){-}C_{(k,1,m{-}1)}(u)\\
\phantom{C_{(k,0,m{+}1)}(u)=}                    {-}C_{(k,0,m{-}1)}(u){-}C_{(k{+}1,0,m{-}1)}(u),\\
C_{(0,l{+}1,m)}(u)=u_2C_{(0,l,m)}(u){-}C_{(1,l{-}1,m{+}2)}(u){-}C_{(1,l{-}2,m{+}2)}(u){-}C_{(0,l{-}1,m)}(u)\\
\phantom{C_{(0,l{+}1,m)}(u)=}                    {-}C_{(2,l{-}1,m)}(u){-}C_{(1,l{+}1,m{-}2)}(u){-}C_{(1,l,m{-}2)}(u),\\
C_{(0,l,m{+}1)}(u)=u_3C_{(0,l,m)}(u){-}C_{(1,l{-}1,m{+}1)}(u){-}C_{(0,l{-}1,m{+}1)}(u){-}C_{(0,l{+}1,m{-}1)}(u)\\
\phantom{C_{(0,l,m{+}1)}(u)=}                    {-}C_{(0,l,m{-}1)}(u){-}C_{(1,l,m{-}1)}(u).\\
\end{gather*}

The lowest $C$-polynomials of $B_3$ were calculated explicitly and arranged in Table~\ref{c_poly_b3}.
\begin{table}[h]
\begin{center}
$\begin{array}{c}
\sharp 0
\\
\hline  \hline
\begin{array}{|l|r|r|r|r|r|r|r|r|r|r|r|r|r|r|}
   C_{(k, l,m)}(u)&1& u_1&u_2&u_1^2& u_3^2&u_1u_2 &u_1u_3^2&u_1^3&u_2^2  &u_2u_3^2&u_1^2u_2  \\ \hline
   C_{(0,0,0)}(u) &1&    &   &     &      &       &       &        &       &       &            \\ \hline
   C_{(1,0,0)}(u) &0&  1 &   &     &      &       &       &        &       &       &                \\ \hline
   C_{(0,1,0)} (u)&0&  0 & 1 &     &      &       &       &        &       &       &             \\ \hline
   C_{(2,0,0)} (u)&-6& 0 &-2 & 1   &      &       &       &        &       &       &        \\ \hline
   C_{(0,0,2)} (u)&-8& -4&-2 & 0   & 1    &       &       &        &       &       &       \\ \hline
   C_{(1,1,0)}(u) &24& 8 & 6 & 0   & -3   & 1     &       &        &       &       &       \\ \hline
   C_{(1,0,2)}(u) &0& -8 &-2 &-4   & 0    & -2    & 1     &        &       &       &           \\ \hline
   C_{(3,0,0)}(u)&-24&-15& -6& 0   & 3    &  -3   & 0     &    1   &       &       &       \\ \hline
   C_{(0,2,0)} (u)&12& 16& 8 & 4   & 0    & 4     &  -2   &    0   &   1   &       &          \\ \hline
   C_{(0,1,2)}(u)&-48&-20&-20& 0   & 6    &  -6   &0      &    0   &  -2   &   1   &          \\ \hline
   C_{(2,1,0)} (u)&0& 8  &-6 & 4   & 0    &  2    & -1    &    0   &  -2   &   0   &    1     \\ \hline
\end{array}
\end{array}$
\\[2ex]
$\begin{array}{c}
\sharp 1
\\
\hline  \hline
\begin{array}{|l|r|r|r|r|r|r|r|}
   C_{(k, l,m)}(u)&u_3& u_1u_3&u_2u_3 & u_3^3&u_1^2u_3 &u_1u_2u_3\\\hline
   C_{(0,0,1)}(u) & 1 &       &       &      &         &        \\\hline
   C_{(1,0,1)}(u) &-3 & 1     &       &      &         &            \\ \hline
   C_{(0,1,1)}(u) & 3 & -2    & 1     &      &         &         \\ \hline
   C_{(0,0,3)}(u) & -9& -3    &-3     & 1    &         &          \\ \hline
   C_{(2,0,1)}(u) & -3&-1     &-2     & 0    &  1      &        \\ \hline
   C_{(1,1,1)}(u)&3 0 & 12    & 8     &  -3  & -2      &    1    \\ \hline
\end{array}
\end{array}$
\end{center}
\caption{Lowest $C$-polynomials of $B_3$ split into two congruence classes $\#=0$ and $\#=1$.}\label{c_poly_b3}
\end{table}

As in the previous case, we can use the Weyl character formula
or generic recurrence relations for $S$-polynomials of $B_3$ valid for $k,l,m\geqslant 2$:
\begin{gather*}
\mathcal{S}_{(k{+}1,l,m)}(u)=u_1\mathcal{S}_{(k,l,m)}(u){-}\mathcal{S}_{(k,l{+}1,m{-}2)}(u){-}\mathcal{S}_{(k,l{-}1,m{+}2)}(u){-}\mathcal{S}_{(k{+}1,l{-}1,m)}(u)\\
\phantom{\mathcal{S}_{(k{+}1,l,m)}(u)=}        {-}\mathcal{S}_{(k{-}1,l{+}1,m)}(u){-}\mathcal{S}_{(k{-}1,l,m)}(u){-}\mathcal{S}_{(k{-}1,l{+}1,m)}(u){-}\mathcal{S}_{(k{-}1,l,m)}(u),\\
\mathcal{S}_{(k,l{+}1,m)}(u)=u_2\mathcal{S}_{(k,l,m)}(u){-}\mathcal{S}_{(k{+}1,l{-}1,m{+}2)}(u){-}\mathcal{S}_{(k{-}1,l{-}1,m{+}2)}(u){-}\mathcal{S}_{(k{+}1,l{-}2,m{+}2)}(u)\\
\phantom{\mathcal{S}_{(k,l{+}1,m)}(u)=}         {-}\mathcal{S}_{(k{+}1,l,m{-}2)}(u){-}\mathcal{S}_{(k{-}2,l{+}1,m)}(u){-}\mathcal{S}_{(k{-}1,l,m{+}2)}(u){-}\mathcal{S}_{(k{-}1,l{+}2,m{-}2)}(u)\\
\phantom{\mathcal{S}_{(k,l{+}1,m)}(u)=}         {-}\mathcal{S}_{(k{+}1,l{+}1,m{-}2)}(u){-}\mathcal{S}_{(k{-}1,l{+}1,m{-}2)}(u){-}\mathcal{S}_{(k{+}2,l{-}1,m)}(u),\\
\mathcal{S}_{(k,l,m{+}1)}(u)=u_3\mathcal{S}_{(k,l,m)}(u){-}\mathcal{S}_{(k{+}1,l{-}1,m{+}1)}(u){-}\mathcal{S}_{(k,l{-}1,m{+}1)}(u){-}\mathcal{S}_{(k{-}1,l,m{+}1)}(u)\\
\phantom{\mathcal{S}_{(k,l,m{+}1)}(u)=}         {-}\mathcal{S}_{(k{-}1,l{+}1,m{-}1)}(u){-}\mathcal{S}_{(k,l{+}1,m{-}1)}(u){-}\mathcal{S}_{(k{+}1,l,m{-}1)}(u){-}\mathcal{S}_{(k,l,m{-}1)}(u).
\end{gather*}

Substitution of variables $x\mapsto u$ transforms vertices of the simplex $F(B_3)$
into vertices of the orthogonality domain $\mathfrak{F}(B_3)=\{P_0,P_1,P_2,P_3\}$ as follows:
\begin{gather*}
(0,0,0)\mapsto (6, 12, 8 )=:P_0,\quad
\omega_1 \mapsto (6, 12, -8 )=:P_1,\\
\tfrac12{\omega_2} \mapsto (-2, -4, 0 )=:P_2,\quad\;\,
\omega_3 \mapsto (-6, 12, 0 )=:P_3.
\end{gather*}
The domain $\mathfrak{F}(B_3)$ and discretization points from Example~\ref{ex_F4_B3} are plotted in Figure~\ref{fig_b3}.

\begin{example}\label{ex_F4_B3}
Let us fix $M=4$. $|F_4(B_3)|=14$ and lattice points $(\frac{s_1}{M},\frac{s_2}{M},\frac{s_3}{M})$ in new coordinates $(u_1,u_2,u_3)$
have the form:
\begin{gather*}
\begin{array}{lll}
(0,0,0)                     \mapsto (6,12,8),       & (0,0,\frac{1}{4})           \mapsto (0,0,2\sqrt2),       & (0,\frac{1}{2},0)           \mapsto (-2,-4,0),\\
(\frac{1}{2},\frac{1}{4},0) \mapsto (2,0,-4),\qquad & (1,0,0)                     \mapsto (6,12,-8),           & (\frac{3}{4},0,0)           \mapsto (4,4,-4\sqrt2),\\
(0,\frac{1}{4},0)           \mapsto (2,0,4),        & (\frac{1}{4},0,0)           \mapsto (4,4,4\sqrt2),       & (\frac{1}{4},\frac{1}{4},0) \mapsto (0,-4,0),\\
(\frac{1}{2},0,0)           \mapsto (2,-4,0),       & (\frac{1}{2},0,\frac{1}{4}) \mapsto (0,0,-2\sqrt2),\qquad& (\frac{1}{4},0,\frac{1}{4}) \mapsto (-2,0,0),\\
(0,\frac{1}{4},\frac{1}{4}) \mapsto (-4,4,0),       & (0,0,\frac{1}{2})           \mapsto (-6,12,0).           &
\end{array}
\end{gather*}
\begin{figure}[h]
\centerline{\includegraphics[scale=0.8]{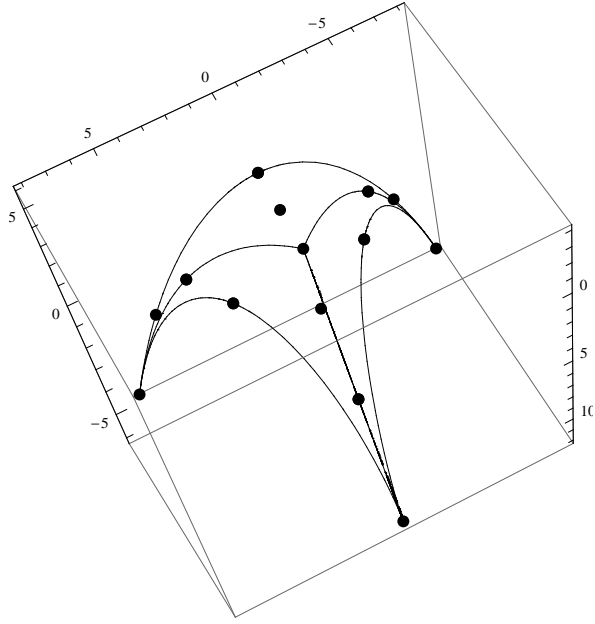}}
\caption{Region of orthogonality $\mathfrak{F}$ and discrete points $\mathfrak{F}_4$ of polynomials of $B_3$.}\label{fig_b3}
\end{figure}
\end{example}

\subsection{$C$- and $S$-polynomials of $C_3$ in three real variables}
Low order $C$-polynomials of group $C_3$ are listed in Table~\ref{c_poly_c3}.
Higher-order polynomials can be obtained from the recurrence relations.
Generic recursions for $k,l,m\geqslant 2$ are:
\begin{gather*}
C_{(k{+}1,l,m)}(u)=u_1C_{(k,l,m)}(u){-}C_{(k,l{-}1,m{+}1)}(u){-}C_{(k{+}1,l{-}1,m)}(u){-}C_{(k{-}1,l{+}1,m)}(u)\\
\phantom{C_{(k{+}1,l,m)}(u)=}        {-}C_{(k{-}1,l,m)}(u){-}C_{(k,l{+}1,m{-}1)}(u),\\
C_{(k,l{+}1,m)}(u)=u_2C_{(k,l,m)}(u){-}C_{(k{+}1,l,m{-}1)}(u){-}C_{(k{-}1,l,m{+}1)}(u){-}C_{(k,l{-}1,m)}(u)\\
\phantom{C_{(k,l{+}1,m)}(u)=}        {-}C_{(k{+}1,l{+}1,m{-}1)}(u){-}C_{(k{-}2,l{+}1,m)}(u){-}C_{(k{-}1,l{-}1,m{+}1)}(u){-}C_{(k{-}1,l{+}1,m{-}1)}(u)\\
\phantom{C_{(k,l{+}1,m)}(u)=}        {-}C_{(k{+}2,l{-}1,m)}(u){-}C_{(k{+}1,l{-}1,m{+}1)}(u){-}C_{(k{+}1,l{-}2,m{+}1)}(u){-}C_{(k{-}1,l{+}2,m{-}1)}(u),\\
C_{(k,l,m{+}1)}(u)=u_3C_{(k,l,m)}(u){-}C_{(k,l{-}2,m{+}1)}(u){-}C_{(k{-}2,l,m{+}1)}(u){-}C_{(k,l{+}2,m{-}1)}(u)\\
\phantom{C_{(k,l,m{+}1)}(u)=}        {-}C_{(k,l,m{-}1)}(u){-}C_{(k{+}2,l{-}2,m{+}1)}(u){-}C_{(k{-}2,l{+}2,m{-}1)}(u){-}C_{(k{+}2,l,m{-}1)}(u).
\end{gather*}

Additional recursions:
\begin{gather*}
C_{(k{+}1,0,0)}(u)=u_1C_{(k,0,0)}(u){-}C_{(k{-}1,1,0)}(u){-}C_{(k{-}1,0,0)}(u),\quad k>1;\\
C_{(k{+}1,l,0)}(u)=u_1C_{(k,l,0)}(u){-}C_{(k,l{-}1,1)}(u){-}C_{(k{-}1,l.0)}(u)\\
\phantom{C_{(k{+}1,l,0)}(u)=}        {-}C_{(k{+}1,l{-}1,0)}(u){-}C_{(k{-}1,l{+}1,0)}(u),\quad k,l>1;\\
C_{(k,l{+}1,0)}(u)=u_2C_{(k,l,0)}(u){-}C_{(k{-}1,l,1)}(u){-}C_{(k,l{-}1,0)}(u){-}C_{(k{+}2,l{-}1,0)}(u){-}C_{(k{-}2,l{+}1,0)}(u)\\
\phantom{C_{(k,l{+}1,0)}(u)=}        {-}C_{(k{-}1,l{-}1,1)}(u){-}C_{(k{+}1,l{-}1,1)}(u){-}C_{(k{+}1,l{-}2,1)}(u),\quad k,l>2;\\
C_{(0,l{+}1,m)}(u)=u_2C_{(0,l,m)}(u){-}C_{(1,l,m{-}1)}(u){-}C_{(0,l{-}1,m)}(u){-}C_{(2,l{-}1,m)}(u){-}C_{(1,l{+}1,m{-}1)}(u)\\
\phantom{C_{(0,l{+}1,m)}(u)=}        {-}C_{(1,l{-}1,m{+}1)}(u){-}C_{(1,l{-}2,m{+}1)}(u),\quad l>2,\ m>1;\\
C_{(0,l,m{+}1)}(u)=u_3C_{(0,l,m)}(u){-}C_{(0,l{-}2,m{+}1)}(u){-}C_{(0,l,m{-}1)}(u){-}C_{(2,l{-}2,m{+}1)}(u)\\
\phantom{C_{(0,l,m{+}1)}(u)=}        {-}C_{(2,l,m{-}1)}(u){-}C_{(0,l{+}2,m{-}1)}(u),\quad l>2,\ m>1;\\
C_{(k{+}1,0,m)}(u)=u_1C_{(k,0,m)}(u){-}C_{(k{-}1,1,m)}(u){-}C_{(k,1,m{-}1)}(u){-}C_{(k{-}1,0,m)}(u),\quad k,m>1;\\
C_{(k,0,m{+}1)}(u)=u_3C_{(k,0,m)}(u){-}C_{(k{-}2,0,m{+}1)}(u){-}C_{(k,0,m{-}1)}(u){-}C_{(k{-}2,2,m{-}1)}(u)\\
\phantom{C_{(k,0,m{+}1)}(u)=}        {-}C_{(k{+}2,0,m{-}1)}(u){-}C_{(k,2,m{-}1)}(u),\quad k>2,\ m>1;\\
C_{(0,l{+}1,0)}(u)=u_2C_{(0,l,0)}(u){-}C_{(0,l{-}1,0)}(u){-}C_{(2,l{-}1,0)}(u){-}C_{(1,l{-}1,1)}(u){-}C_{(1,l{-}2,1)}(u),\;\; l>2;\\
C_{(0,0,m{+}1)}(u)=u_3C_{(0,0,m)}(u){-}C_{(0,0,m{-}1)}(u){-}C_{(2,0,m{-}1)}(u){-}C_{(0,2,m{-}1)}(u),\quad m>1.
\end{gather*}

Generic recursions for $S$-polynomials hold true when $k,l,m\geqslant 2$
\begin{gather*}
\mathcal{S}_{(k{+}1,l,m)}(u)=u_1\mathcal{S}_{(k,l,m)}(u){-}\mathcal{S}_{(k,l{-}1,m{+}1)}(u){-}\mathcal{S}_{(k{+}1,l{-}1,m)}(u){-}\mathcal{S}_{(k{-}1,l{+}1,m)}(u)\\
\phantom{\mathcal{S}_{(k{+}1,l,m)}(u)=}        {-}\mathcal{S}_{(k{-}1,l,m)}(u){-}\mathcal{S}_{(k,l{+}1,m{-}1)}(u),\\
\mathcal{S}_{(k,l{+}1,m)}(u)=u_2\mathcal{S}_{(k,l,m)}(u){-}\mathcal{S}_{(k{+}1,l,m{-}1)}(u){-}\mathcal{S}_{(k{-}1,l,m{+}1)}(u){-}\mathcal{S}_{(k{+}2,l{-}1,m)}(u)\\
\phantom{\mathcal{S}_{(k,l{+}1,m)}(u)=}        {-}\mathcal{S}_{(k{-}2,l{+}1,m)}(u){-}\mathcal{S}_{(k{-}1,l{-}1,m{+}1)}(u){-}\mathcal{S}_{(k{-}1,l{+}1,m{-}1)}(u){-}\mathcal{S}_{(k{+}1,l{+}1,m{-}1)}(u)\\
\phantom{\mathcal{S}_{(k,l{+}1,m)}(u)=}        {-}\mathcal{S}_{(k{+}1,l{-}1,m{+}1)}(u){-}\mathcal{S}_{(k{+}1,l{-}2,m{+}1)}(u){-}\mathcal{S}_{(k{-}1,l{+}2,m{-}1)}(u){-}\mathcal{S}_{(k,l{-}1,m)}(u),\\
\mathcal{S}_{(k,l,m{+}1)}(u)=u_3\mathcal{S}_{(k,l,m)}(u){-}\mathcal{S}_{(k,l{-}2,m{+}1)}(u){-}\mathcal{S}_{(k{-}2,l,m{+}1)}(u){-}\mathcal{S}_{(k,l{+}2,m{-}1)}(u)\\
\phantom{\mathcal{S}_{(k,l,m{+}1)}(u)=}        {-}\mathcal{S}_{(k{+}2,l{-}2,m{+}1)}(u){-}\mathcal{S}_{(k{-}2,l{+}2,m{-}1)}(u){-}\mathcal{S}_{(k{+}2,l,m{-}1)}(u){-}\mathcal{S}_{(k,l,m{-}1)}(u).
\end{gather*}

\begin{table}[h]
\begin{center}
$\begin{array}{c}
\sharp 0     \\
\hline  \hline
\begin{array}{|l|r|r|r|r|r|r|r|r|r|r|r|r|r|r|}
    C_{(k, l,m)}(u) &1 &u_2&u_1^2&u_1u_3& u_2^2&u_1^2u_2&u_3^2&u_1u_2u_3 &u_2^3  &u_2u_3^2    \\\hline
    C_{(0,0,0)} (u)   &1 &   &     &      &      &        &     &          &       &            \\\hline
    C_{(0,1,0)} (u)   &0 & 1 &     &      &      &        &     &          &       &          \\\hline
    C_{(2,0,0)} (u)   &-6&-2 &1    &      &      &        &     &          &       &         \\\hline
    C_{(1,0,1)} (u)   & 0&-2 &0    &  1   &      &        &     &          &       &            \\ \hline
    C_{(0,2,0)} (u)   &12& 8 & -4  &  -2  &   1  &        &     &          &       &            \\ \hline
    C_{(2,1,0)} (u)   & 0&-6 &  0  &   -1 &   -2 &    1   &     &          &       &            \\\hline
    C_{(0,0,2)} (u)   &-8& -8&  4  &   4  &   -2 &   0    & 1   &          &       &            \\\hline
    C_{(1,1,1)} (u)   & 0& 12& 0   &   -4 &   4  &   -2   &  -3 & 1        &       &       \\\hline
    C_{(0,3,0)} (u)   & 0& 9 & 0   &   3  &   6  &   -3   &  3  &-3        &1      &                     \\\hline
    C_{(0,1,2)} (u)   & 0&-18&  0  &  3   &   -12&    6   & 3   & 3        &-2     &1    \\\hline
\end{array}
\end{array}$
\\[3ex]
$\begin{array}{c}
\sharp 1    \\
\hline  \hline
\begin{array}{|l|r|r|r|r|r|r|r|r|r|r|r|r|r|rrrrrrrrrrrrrrrr}
       C_{(k, l,m)}(u) &u_1  &u_3& u_1u_2& u_1^3&u_2u_3&u_1^2u_3 &u_1u_2^2   &u_1u_3^2   &u_2^2u_3   &v_3^3\\\hline
       C_{(1,0,0)}   (u) & 1   &   &       &      &      &          &          &           &           &\\\hline
       C_{(0,0,1)}(u)    & 0   & 1 &       &      &      &          &          &           &           &      \\\hline
       C_{(1,1,0)}(u)    &-4   &-3 &   1   &      &      &          &          &           &           &\\\hline
       C_{(3,0,0)} (u)   &-3   &3  &   -3  &  1   &     &           &          &           &           &   \\ \hline
       C_{(0,1,1)} (u)   & 4   & 6 &   -2  &  0   &  1  &           &          &           &           &   \\ \hline
       C_{(2,0,1)}(u)    &  0  & -9&   0   &   0  &  -2 &  1        &          &           &           &   \\\hline
       C_{(1,2,0)}(u)    & 12  & -3&   9   &   -4 & -1  &   -2      & 1        &           &           &     \\\hline
       C_{(1,0,2)}(u)    &-12  & -6&   -6  &  4   &  -1 &  4        &-2        &1          &           &          \\\hline
       C_{(0,2,1)}(u)    & 0   & 27&    0  &  0   &  12 & -6        & 0        &-2         &1          &               \\\hline
       C_{(0,0,3)} (u)   & 0   &-27&    0  &  0   & -18 &  9        &0         & 6         &-3         &1    \\\hline
\end{array}
\end{array}$
\end{center}
\caption{Lowest $C$-polynomials of $C_3$ split into two congruence classes $\#=0$ and $\#=1$.}\label{c_poly_c3}
\end{table}

Let us study the transformation of the fundamental region $F(C_3)\to \mathfrak{F}(C_3)$.
During the substitution $x\mapsto u$, the vertices of the simplex $F(C_3)$ (see Appendix) go to the points $\{P_0,P_1,P_2,P_3\}$:
\begin{gather*}
(0,0,0)  \mapsto  (6, 12, 8 ) =:P_0,\quad
\omega_1 \mapsto  (2, -4, -8 )=:P_1,\\
\omega_2 \mapsto  (-2, -4, 8 )=:P_2,\qquad
\omega_3 \mapsto  (-6, 12, -8 )=:P_3.
\end{gather*}
The region of orthogonality of polynomials of $C_3$ is presented in Figure~\ref{fig_c3}.

\begin{example}
Let us fix $M=4$. $|F_4(C_3)|=14$ and lattice points $(\frac{s_1}{M},\frac{s_2}{M},\frac{s_3}{M})\in F_4$
map to $(u_1,u_2,u_3)$:
\begin{gather*}
\begin{array}{lll}
(0,0,0)                     \mapsto (6,12,8),                  & (0,0,\frac{1}{4})           \mapsto (3\sqrt2,6,2\sqrt2),       & (0,\frac{1}{4},0)           \mapsto (2,0,0),\\
(0,\frac{1}{4},\frac{1}{4}) \mapsto (-\sqrt2,-2,2\sqrt2),\quad & (\frac{1}{4},0,\frac{1}{4}) \mapsto (\sqrt2,-2,-2\sqrt2), \quad& (\frac{1}{4},0,\frac{1}{2}) \mapsto (-2,0,0),\\
(\frac{1}{4},\frac{1}{4},0) \mapsto (0,-4,0),                  & (0,\frac{1}{2},0)           \mapsto (-2,-4,8),                 & (0,\frac{1}{4},\frac{1}{2}) \mapsto (-4,4,0),\\
(\frac{1}{2},0,0)           \mapsto (2,-4,-8),                 & (0,0,1)                     \mapsto (-6,12,-8),                & (0,0,\frac{1}{2})           \mapsto (0,0,0),\\
(\frac{1}{4},0,0)           \mapsto (4,4,0),                   & (0,0,\frac{3}{4})           \mapsto (-3\sqrt2,6,-2\sqrt2).     &
\end{array}
\end{gather*}
These points are shown on the orthogonality domain $\mathfrak{F}(C_3)$ in Figure~\ref{fig_c3}.

\begin{figure}[h]
\centerline{\includegraphics[scale=0.8]{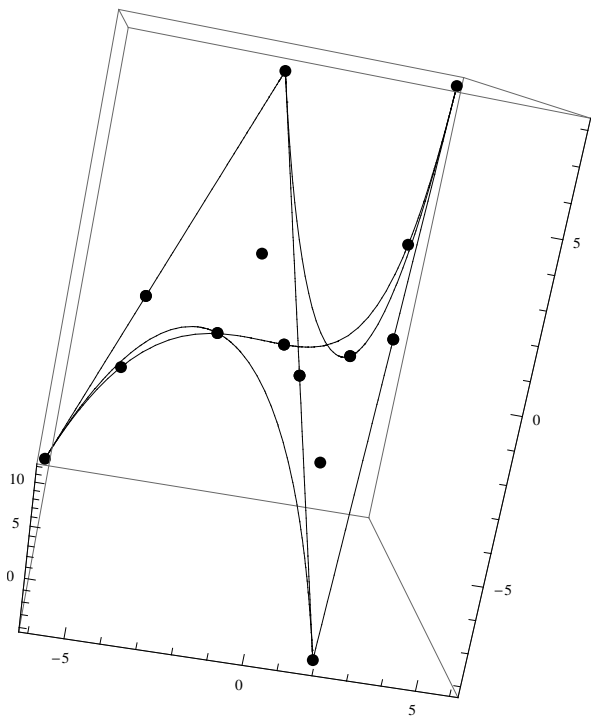}}
\caption{Region of orthogonality $\mathfrak{F}$ and discrete points $\mathfrak{F}_4$ of polynomials of $C_3$.}\label{fig_c3}
\end{figure}
\end{example}

\section{Branching rules for polynomials of $A_2$, $C_2$, $G_2$, $A_3$, $B_3$, and $C_3$.}\label{sec_reductions}
Consider simple Lie groups of rank $\leq3$.
An uncommon transformation (`branching') of polynomials of groups into the sum of polynomials of their subgroups.
We are interested in cases when a complete reduction of all polynomials of the group is achieved,
namely when the sum consists of complete polynomials of the subgroup.

Probably the most interesting class of cases can be identified by the maximal semisimple subgroups contained in the simple Lie groups. There are the following cases to consider:
\begin{gather*}
A_2\supset A_1;
\\
C_2\supset A_1\times A_1, \quad
C_2\supset A_1;
\\
G_2\supset A_1\times A_1, \quad
G_2\supset A_2, \quad
G_2\supset A_1.
\\[2ex]
A_3\supset C_2, \quad
A_3\supset A_1\times A_1;
\\
B_3\supset A_3, \quad
B_3\supset A_1\times  A_1\times A_1, \quad
B_3\supset G_2;
\\
C_3\supset C_2\times A_1, \quad
C_3\supset A_2, \quad
C_3\supset A_1.
\end{gather*}

In addition, complete reduction of polynomials takes place when the subgroup is non-semisimple reductive.
Three such subgroups are maximal with rank $\leq3$,
namely $A_2\supset A_1\times U_1$, $A_3\supset A_2\times U_1$, and $B_3\supset B_2\times U_1$.
Here, irreducible polynomials of the 1-parametric group $U_1$ are trivially all monomials, say $Y^k$, with integer $k$ positive or negative.

Finally, complete reduction of polynomials happens when the relation between the pair of groups is not inclusion but subjoining \cite{PSS}.
A classification of all such relations is in~\cite{MPi}.

\subsection{Polynomials of the Lie algebra $A_1$}
It was shown in~\cite{NPT0, NPT} that there is a one-to-one correspondence between $C$-
and $S$-polynomials of $A_1$ and Chebyshev polynomials of the first and second kind respectively.
Here we just write down the explicit form of the $C$-polynomials of $A_1$ (the Chebyshev polynomials of the first kind)
traditionally denoted by $T_m$, $m=0,1,2,\dots$.

\begin{gather*}
T_0(x)=1,\quad
T_1(x)=x,\quad
T_2(x)=2x^2-1,\quad
T_3(x)=4x^3-3x,\quad
T_4(x)=8x^4-8x^2+1,\\
T_{m+1}(x)=2xT_m-T_{m-1}.\label{def_T0_T1_Tn}
\end{gather*}

The equivalent set of polynomials can be obtained via the substitution $x=\frac y 2$, and up to the scaling factor they are:
\begin{gather*}
\begin{gathered}
\tilde{T}_0(y)=1,\quad
\tilde{T}_1(y)=y,\quad
\tilde{T}_2(y)=y^2-1,\quad
\tilde{T}_3(y)=y^3-3y,\quad
\tilde{T}_4(y)=y^4-4y^2+2, \dots
\end{gathered}
\end{gather*}

\subsection{Reductions of polynomials of rank two simple Lie groups.}
In this section, we present explicit forms of projection matrices ${\rm Pr}$
providing reductions from the rank two simple Lie groups $G$ to their maximal subgroups $H$.
For the each set $\{G\supset H,{\rm Pr}\}$ we obtain the corresponding substitution of variables
$X_i\mapsto f_i(Y_1,Y_2)$ or $X_i\mapsto f_i(Y)$, $i=1,2$
that reduces the orthogonal polynomials of $G$ to the polynomials of $H$.
Explicit examples of $C$-polynomial reductions are shown for all pairs $G\supset H$.

Here we use the notation of polynomial variables taken from the paper~\cite{NPT},
namely \mbox{$X_i:=C_{\w_i}$}, instead of $u_i=X_i$, $i=1,2$ as in previous sections.

\subsection*{Reduction from $A_2$ to $A_1$.}
\begin{gather*}
{\rm Pr}=\left(
\begin{array}{cc}
2&2
\end{array}
\right),
\\[2ex]
\begin{array}{@{}l}
C_{(1,0)}(A_2) \ \mapsto \  C_{(2)}(A_1){+}C_{(0)}(A_1),\\[1ex]
C_{(0,1)}(A_2) \ \mapsto \  C_{(2)}(A_1){+}C_{(0)}(A_1),
\end{array}
\qquad \text{therefore} \qquad
\begin{array}{l}
X_1 \ \mapsto \ Y^2{-}1,\\
X_2 \ \mapsto \ Y^2{-}1.
\end{array}
\end{gather*}

\noindent{Example of reduction: }
Consider the reduction of a few low order polynomials of $A_2$.
Hereafter the explicit forms of polynomials of simple Lie groups of rank two are taken from the paper~\cite{NPT}.
\begin{align*}
C_{(1,1)}(A_2)=&X_1X_2{-}3 \ \mapsto \ Y^4{-}2Y^2{-}2=\tilde{T}_4(Y){+}2\tilde{T}_2(Y),\\
C_{(0,2)}(A_2)=&X_2^2{-}2X_1 \ \mapsto \ Y^4{-}4Y^2{+}3=\tilde{T}_4(Y){+}1.
\end{align*}

\subsection*{Reduction from $C_2$ to $A_1\times A_1=\colon A_1^{(1)}\times A_1^{(2)}$.}
\begin{gather*}
{\rm Pr}=\left(
\begin{array}{cc}
1&1\\
0&1
\end{array}
\right),
\\[2ex]
\begin{array}{@{}l}
C_{(1,0)}(C_2) \ \mapsto \  C_{(1)}(A_1^{(1)})C_{(0)}(A_1^{(2)})+C_{(0)}(A_1^{(1)})C_{(1)}(A_1^{(2)}),\\[1ex]
C_{(0,1)}(C_2) \ \mapsto \  C_{(1)}(A_1^{(1)})C_{(1)}(A_1^{(2)}),
\end{array}
\qquad \text{therefore} \qquad
\begin{array}{l}
X_1 \ \mapsto \  Y_1+Y_2,\\
X_2 \ \mapsto \  Y_1 Y_2.
\end{array}
\end{gather*}

{\samepage
\noindent{Example of reduction: }
In this example, Chebyshev polynomials $\tilde{T}_m(Y_1)$, $m=1,2,\dots$ correspond to the Lie group $A_1^{(1)}$,
and $\tilde{T}_m(Y_2)$ correspond to group $A_1^{(2)}$.
\begin{align*}
C_{(2,1)}(C_2){=}&X_1^2X_2{-}2X_2^2{-}6X_2 \ \mapsto \  Y_1^3Y_2{+}Y_1Y_2^3{-}6Y_1Y_2{=} \tilde{T}_3(Y_1)Y_2{+}Y_1\tilde{T}_3(Y_2),\\
C_{(3,0)}(C_2){=}&X_1^3{-}3X_1X_2{-}3X_1 \ \mapsto \  Y_1^3{+}Y_2^3{-}3Y_1{-}3Y_2=\tilde{T}_3(Y_1){+}\tilde{T}_3(Y_2).
\end{align*}
}

\subsection*{Reduction from $C_2$ to $A_1$.}
\begin{gather*}
{\rm Pr}=\left(
\begin{array}{cc}
3&4
\end{array}
\right),
\\[2ex]
\begin{array}{@{}l}
C_{(1,0)}(C_2) \ \mapsto \  C_{(3)}(A_1){+}C_{(1)}(A_1),\\
C_{(0,1)}(C_2) \ \mapsto \  C_{(4)}(A_1){+}C_{(2)}(A_1),
\end{array}
\qquad \text{therefore} \qquad
\begin{array}{l}
X_1 \  \mapsto \  Y^3{-}2Y,\\
X_2 \ \mapsto \  Y^4{-}3Y^2.
\end{array}
\end{gather*}

\noindent{Example of reduction: }
\begin{align*}
C_{(1,1)}(C_2){=}&X_1X_2{-}2X_1 \ \mapsto \  Y^7{-}5Y^5{+}4Y^3{+}4Y{=}\tilde{T}_7(Y){-}2\tilde{T}_5(Y){+}21Y,\\
C_{(2,0)}(C_2){=}&X_1^2{-}2X_2{-}4 \ \mapsto  \  Y^6{-}6Y^4{+}10Y^2{-}4{=}\tilde{T}_6(Y){+}\tilde{T}_2(Y).
\end{align*}

\subsection*{Reduction from $G_2$ to $A_2$.}
\begin{gather*}
{\rm Pr}=\left(
\begin{array}{cc}
1&0\\
1&1
\end{array}
\right),
\\[2ex]
\begin{array}{@{}l}
C_{(1,0)}(G_2) \ \mapsto \  C_{(1,1)}(A_2),\\
C_{(0,1)}(G_2) \ \mapsto \  C_{(1,0)}(A_2){+}C_{(0,1)}(A_2),
\end{array}
\qquad \text{therefore} \qquad
\begin{array}{l}
X_1 \ \mapsto \   Y_1Y_2{-}3,\\
X_2 \ \mapsto \   Y_1{+}Y_2.
\end{array}
\end{gather*}

\noindent{Example of reduction: }
Below $C_{(a,b)}(\bar{Y}):=C_{(a,b)}(Y_1,Y_2)$ are the polynomials of $A_2$.
\begin{align*}
C_{(0,2)}(G_2)=&X_2^2{-}2X_2{-}2X_1{-}6 \ \mapsto \
Y_1^2{+}Y_2^2{-}2Y_1{-}2Y_2=C_{(2,0)}(\bar Y){+}C_{(0,2)}(\bar Y),\\
C_{(1,1)}(G_2)=&-2X_2^2{+}X_1X_2{+}2X_2{+}4X_1{+}12 \ \mapsto \
Y_1^2Y_2{+}Y_1Y_2^2{-}2Y_1^2{-}2Y_2^2{-}Y_1{-}Y_1=\\=&C_{(2,1)}(\bar Y){+}C_{(1,2)}(\bar Y).
\end{align*}

\subsection*{Reduction from $G_2$ to $A_1\times A_1=\colon A_1^{(1)}\times A_1^{(2)}$.}
\begin{gather*}{\rm Pr}=\left(
\begin{array}{cc}
1&1\\
3&1
\end{array}
\right),
\\[2ex]
\begin{array}{@{}l}
C_{(1,0)}(G_2) \ \mapsto \  C_{(1)}(A_1^{(1)})C_{(3)}(A_1^{(2)})+ C_{(2)}(A_1^{(1)})C_{(0)} (A_1^{(2)}),\\[1ex]
C_{(0,1)}(G_2) \ \mapsto \  C_{(1)}(A_1^{(1)})C_{(1)}(A_1^{(2)})+C_{(0)}(A_1^{(1)})C_{(2)}(A_1^{(2)}),\end{array}\\
\qquad \text{therefore} \qquad
\begin{array}{l}
X_1 \ \mapsto \  Y_1 Y_2^3{-}3Y_1Y_2{+}Y_1^2{-}2,\\
X_2 \ \mapsto \  Y_1Y_2{+}Y_2^2{-}2.
\end{array}
\end{gather*}

\noindent{Example of reduction: }
Here each polynomial $\tilde{T}_m(Y_1)$, $m=1,2,\dots$ corresponds to $A_1^{(1)}$
and $\tilde{T}_m(Y_2)$ does to $A_1^{(2)}$.
\begin{align*}
C_{(0,2)}(G_2)=&X_2^2{-}2X_2{-}2X_1{-}6 \ \mapsto \
Y_2^4{+}Y_1^2Y_2^2{-}6Y_2^2{-}2Y_1^2{+}6=\\
=&\tilde{T}_4(Y_2){+}\tilde{T}_2(Y_1)\tilde{T}_2(Y_2){+}2,
\end{align*}
\begin{align*}
C_{(1,1)}(G_2)=&-2X_2^2{+}X_1X_2{+}2X_2{+}4X_1{+}12 \ \mapsto \ \\
\mapsto& \  Y_1Y_2^5{+}Y_1^2Y_2^4{-}2Y_2^4{-}5Y_1Y_2^3{-}4Y_1^2Y_2^2{+}8Y_2^2{+}Y_1^3Y_2{+}2Y_1Y_2{+}Y_1^2{-}4=\\
=&Y_1\tilde{T}_5(Y_2){-}\tilde{T}_3(Y_1)Y_2{+}\tilde{T}_2(Y_1)\tilde{T}_4(Y_2).
\end{align*}

\subsection*{Reduction from $G_2$ to $A_1$.}
\begin{gather*}
{\rm Pr}=\left(
\begin{array}{cc}
10&6
\end{array}
\right),
\\[2ex]
\begin{array}{@{}l}
C_{(1,0)}(G_2) \ \mapsto \  C_{(10)}(A_1){+}C_{(8)}(A_1){+}C_{(2)}(A_1),\\
C_{(0,1)}(G_2) \ \mapsto \  C_{(6)}(A_1){+}C_{(4)}(A_1){+}C_{(2)}(A_1),
\end{array}\\
\qquad \text{therefore} \qquad
\begin{array}{l}
X_1 \ \mapsto \  Y^{10}{-}9Y^8{+}27Y^6{-}31Y^4{+}10Y^2{-}2,\\
X_2 \ \mapsto \  Y^6{-}5Y^4{+}6Y^2{-}2.
\end{array}
\end{gather*}

\noindent{Example of reduction: }
Because the polynomials of $G_2$ and the variable transformations are too cumbersome in this case
we present only one example of polynomial reduction.
\begin{align*}
C_{(0,2)}(G_2)=&X_2^2{-}2X_2{-}2X_1{-}6 \ \mapsto \
Y^{12}{-}12 Y^{10}{+}55Y^8{-}120Y^6{+}128Y^4{-}56Y^2{+}6{=}\\
{=}&\tilde{T}_{12}(Y){+}\tilde{T}_8(Y){+}3\tilde{T}_4(Y){+}8\tilde{T}_2(Y){+}12.
\end{align*}

\subsection{Reductions of polynomials of rank three simple Lie groups.}
In this section, we obtain the projection matrices ${\rm Pr}$
providing reductions from the rank three simple Lie groups $G$ to their maximal subgroups $H$.
For each set $\{G\supset H,{\rm Pr}\}$ we also obtain the corresponding substitution of variables
$X_i\mapsto f_i(Y_1,Y_2,Y_3)$ or $X_i\mapsto f_i(Y_1,Y_2)$, or $X_i\mapsto f_i(Y)$, $i=1,2,3$.
This reduces the orthogonal polynomials of $G$ to the polynomials~of~$H$.
Explicit examples of $C$-polynomial reductions are shown for all pairs $G\supset H$.

Note that in this section we use the notation of polynomial variables taken from the paper~\cite{NPT},
namely $X_i:=C_{\w_i}$, instead of $u_i=X_i$, $i=1,2,3$ as in the beginning of the~paper.

\subsection*{Reduction from $A_3$ to $C_2$.}
\begin{gather*}
{\rm Pr}=\left(
\begin{array}{ccc}
1&0&1\\
0&1&0
\end{array}\right);
\\[2ex]
\begin{array}{@{}l}
      C_{(1,0,0)}(A_3)\mapsto C_{(1,0)}(C_2),\\
      C_{(0,1,0)}(A_3)\mapsto C_{(0,1)}(C_2){+}2C_{(0,0)}(C_2),\\
      C_{(0,0,1)}(A_3)\mapsto C_{(1,0)}(C_2);
    \end{array} \qquad \text{therefore} \qquad \begin{array}{l}X_1\mapsto\ Y_1,\\
  X_2\mapsto\ Y_2{+}2,\\
  X_3\mapsto\ Y_1.\end{array}\end{gather*}

\noindent{Example of reduction: }
Below $C_{(a,b)}(\overline{Y}):=C_{(a,b)}(Y_1,Y_2)$  are the polynomials of the subalgebra~$C_2$.
\begin{gather*}
  C_{(1,0,1)}(A_3) = -4 + X_1X_3 \mapsto -4+Y_1^2=C_{(2,0)}(\overline{Y})+2 Y_2;\\
  C_{(0,1,1)}(A_3) =-3X_1+X_2X_3 \mapsto -Y_1+Y_1Y_2 = C_{(1,1)}(\overline{Y}) + Y_1.\\
\end{gather*}

\subsection*{Reduction from $A_3$ to $A_1 \times A_1=:A_1^{(1)} \times A_1^{(2)}$.}
\begin{gather*}
{\rm Pr}=\left(
\begin{array}{ccc}
1&0&1\\
1&2&1
\end{array}\right);
\\[2ex]
\begin{array}{@{}l}
C_{(1,0,0)}(A_3)\mapsto C_{1}(A_1^{(1)})C_{1}(A_1^{(2)}),\\
C_{(0,1,0)}(A_3)\mapsto C_0(A_1^{(1)})C_{2}(A_1^{(2)})+2C_0(A_1^{(1)})C_{0}(A_1^{(2)}),\\
C_{(0,0,1)}(A_3)\mapsto C_{1}(A_1^{(1)})C_{1}(A_1^{(2)});
\end{array} \qquad{\rm therefore} \qquad \begin{array}{l}
 X_1\mapsto\ Y_1Y_2,\\
  X_2\mapsto\ Y_2^2,\\
  X_3\mapsto\ Y_1Y_2.\end{array}\end{gather*}

\noindent{Example of reduction: }
In this example, Chebyshev polynomials $\tilde{T}_m(Y_i)$, $m=1,2,\dots$ correspond to the Lie group $A_1^{(i)}, i=1,2.$
\begin{gather*}
  C_{(1,0,1)} = -4 + X_1X_3 \mapsto -4+Y_1^2Y_2^2=\tilde{T}_2(Y_1)\tilde{T}_2(Y_2)+2\tilde{T}_2(Y_1)+2\tilde{T}_2(Y_2);\\
  C_{(0,1,1)} =-3X_1+X_2X_3 \mapsto -3Y_1Y_2+Y_1Y_2^3 = Y_1\tilde{T}_3(Y_2).
\end{gather*}

\subsection*{Reduction from $B_3$ to $A_3$.}
\begin{gather*}
{\rm Pr} = \left(
\begin{array}{ccc}
0&1&1\\
1&1&0\\
0&1&0
\end{array}\right);
\\[2ex]
\begin{array}{@{}l}
C_{(1,0,0)}(B_3)\mapsto C_{(0,1,0)}(A_3){+}C_{(1,0,1)}(A_3),\\
      C_{(0,1,0)}(B_3)\mapsto C_{(1,1,1)}(A_3)+C_{(0,1,2)}(A_3),\\
      C_{(0,0,1)}(B_3)\mapsto C_{(1,0,0)}(A_3)+C_{(0,0,1)}(A_3)+C_{(0,1,1)}(A_3);
      \end{array} \\ \qquad \text{therefore} \qquad \begin{array}{l}
X_1\mapsto\ Y_2+Y_1Y_3-4,\\
X_2\mapsto\ Y_1Y_2Y_3+Y_2Y_3^2-3Y_1^2-3Y_3^2-2Y_2^2-Y_1Y_3+4Y_2+4,\\
X_3\mapsto\ Y_2Y_3-2Y_1+Y_3.\end{array}\end{gather*}

\noindent{Example of reduction: }
Below $C_{(a,b,c)}(\overline{Y}):=C_{(a,b,c)}(Y_1,Y_2,Y_3)$  are polynomials of the subalgebra $A_3$.
\begin{gather*}
C_{(1,0,1)}(B_3) = X_1X_3 - 3 X_3 \mapsto
   Y_1Y_2Y_3^2-2Y_1^2Y_3+Y_1Y_3^2+Y_2^2Y_3-2Y_1Y_2-6Y_2Y_3+14Y_1-7Y_3 = \\ =C_{(1,1,2)}(\overline{Y})+4C_{(1,0,2)}(\overline{Y})+C_{(0,2,1)}(\overline{Y})+2C_{(2,0,1)}(\overline{Y})-7C_{(0,1,1)}(\overline{Y})+5C_{(1,1,0)}(\overline{Y})-9Y_1+6Y_3;\\
C_{(2,0,0)}(B_3) = X_1^2-2X_2-6 \mapsto  Y_1^2Y_3^2-2Y_2Y_3^2+6Y_1^2+5Y_2^2+6Y_3^2-6Y_1Y_3-16Y_2+2=\\
=C_{(2,0,2)}(\overline{Y})+2C_{(2,1,0)}(\overline{Y})+5C_{(0,2,0)}(\overline{Y})+6C_{(2,0,0)}(\overline{Y})+6C_{(0,0,2)}(\overline{Y})+9C_{(1,0,1)}(\overline{Y})+8Y_2+24.
\end{gather*}

\subsection*{Reduction from $B_3$ to $3A_1=:A_1^{(1)}\times A_1^{(2)}\times A_1^{(3)}$.}
\begin{gather*}{\rm Pr}=\left(
\begin{array}{ccc}
1&1&0\\
1&1&1\\
0&2&1
\end{array}\right);
\\[2ex]
\begin{array}{@{}l}
      C_{(1,0,0)}(B_3)\mapsto  C_1(A_1^{1})C_1(A_1^{2})C_0(A_1^{3})+C_0(A_1^{1})C_0(A_1^{2})C_2(A_1^{3}), \\
      C_{(0,1,0)}(B_3)\mapsto  C_1(A_1^{1})C_1(A_1^{2})C_2(A_1^{3})+C_0(A_1^{1})C_2(A_1^{2})C_0(A_1^{3})+C_2(A_1^{1})C_0(A_1^{2})C_0(A_1^{3}), \\
      C_{(0,0,1)}(B_3)\mapsto  C_0(A_1^{1})C_1(A_1^{2})C_1(A_1^{3})+C_1(A_1^{1})C_0(A_1^{2})C_1(A_1^{3});
\end{array}\qquad \\ \text{therefore} \qquad \begin{array}{l}
  X_1\mapsto\ Y_1Y_2+Y_3^2-2,\\
  X_2\mapsto\ Y_1Y_2(Y_3^2-2)+Y_1^2+Y_2^2-4,\\
  X_3\mapsto\ Y_1Y_3+Y_2Y_3.\end{array}\end{gather*}

\noindent{Example of reduction: }
The variables refer to the three different subalgebras $A_1$.
The Chebyshev polynomials $\tilde{T}_m(Y_i)$, $m=1,2,\dots$ correspond to the Lie group $A_1^{(i)}, i=1,2,3.$
\begin{gather*}
C_{(1,0,1)}(B_3) = X_1X_3 - 3 X_3 \mapsto Y_1^2Y_2Y_3+Y_1Y_2^2Y_3+Y_1Y_3^3+Y_2Y_3^3-5Y_1Y_3-5Y_2Y_3= \\ =Y_1\tilde{T}_3(Y_3)+Y_2\tilde{T}_3(Y_3)+\tilde{T}_2(Y_1)Y_2Y_3+Y_1\tilde{T}_2(Y_2)Y_3;\\
C_{(0,1,1)}(B_3) = X_2X_3-2X_1X_3+3X_3 \mapsto \\ Y_1^2Y_2Y_3^3+Y_1Y_2^2Y_3^3-3Y_1^2Y_2Y_3-3Y_1Y_2^2Y_3+Y_1^3Y_3+Y_2^3Y_3-2Y_1Y_3^3-2Y_2Y_3^3+3Y_1Y_3+3Y_2Y_3=\\=\tilde{T}_2(Y_1)Y_2\tilde{T}_3(Y_3)+Y_1\tilde{T}_2(Y_2)\tilde{T}_3(Y_3)+\tilde{T}_3(Y_1)Y_3
+\tilde{T}_3(Y_2)Y_3.
\end{gather*}

\subsection*{Reduction from $B_3$ to $G_2$.}
\begin{gather*}{\rm Pr}=\left(
\begin{array}{ccc}
0&1&0\\
1&0&1
\end{array}
\right);
\\[2ex]
\begin{array}{@{}l}
C_{(1,0,0)}(B_3)\mapsto C_{(0,1)}(G_2),\\
C_{(0,1,0)}(B_3)\mapsto C_{(1,0)}(G_2)+C_{(0,1)}(G_2), \\
C_{(0,0,1)}(B_3)\mapsto C_{(0,1)}(G_2)+C_{(0,0)}(G_2);
\end{array}\qquad \text{therefore} \qquad
\begin{array}{l}
  X_1\mapsto\ Y_2,\\
  X_2\mapsto\ Y_1+Y_2,\\
  X_3\mapsto\ Y_2+2.
\end{array}\end{gather*}

\noindent{Example of reduction: }
In this example, $C_{(a,b,c)}(\overline{Y}):=C_{(a,b)}(Y_1,Y_2)$  are the polynomials of subalgebra~$G_2$.

\begin{gather*}
  C_{(1,0,1)}(B_3) = X_1X_3 - 3 X_3 \mapsto Y_2^2 - Y_2-6 = C_{(0,2)}(\overline{Y})+2Y_1+Y_2;\\
  C_{(1,1,0)}(B_3) =X_1^2 -2X_2-6 \mapsto Y_2^2-2Y_1-2Y_2-6 = C_{(0,2)}(\overline{Y}).
\end{gather*}

\subsection*{Reduction from $C_3$ to $C_2 \times A_1$.}
\begin{gather*}{\rm Pr}=\left(
\begin{array}{ccc}1&0&0\\
0&1&1\\
0&0&1
\end{array}\right);
\\[2ex]
\begin{array}{@{}l}
C_{(1,0,0)}(C_3)\mapsto C_{(1,0)}(C_2)C_0(A_1)+C_{(0,0)}(C_2)C_1(A_1),\\
      C_{(0,1,0)}(C_3)\mapsto C_{(0,1)}(C_2)C_0(A_1)+C_{(1,0)}(C_2)C_1(A_1),\\
      C_{(0,0,1)}(C_3)\mapsto C_{(0,1)}(C_2)C_1(A_1);
\end{array}\qquad \text{therefore} \qquad \begin{array}{l}
  X_1\mapsto\ Y_1+Y_3,\\
  X_2\mapsto\ Y_2+Y_1Y_3,\\
  X_3\mapsto\ Y_2Y_3.\end{array}\end{gather*}

\noindent{Example of reduction: }
Let $C_{(a,b,c)}(\overline{Y}):=C_{(a,b)}(Y_1,Y_2)$  be the polynomials of the subalgebra $C_2,$ and Chebyshev polynomials $\tilde{T}_m(Y_3)$, $m=1,2,\dots$ correspond to the Lie group $A_1.$
\begin{gather*}
  C_{(1,0,1)}(C_3) = X_1X_3-2X_2 \mapsto Y_1Y_2Y_3+Y_2Y_3^2-2Y_1Y_3-2Y_2=C_{(1,1)}(\overline{Y})Y_3+Y_2\tilde{T}_2(Y_3);\\
  C_{(2,0,0)}(C_3) = X_1^2-2X_2-6 \mapsto Y_1^2+Y_3^2-2Y_2-6=C_{(2,0)}(\overline{Y})+\tilde{T}_2(Y_3).
\end{gather*}

\subsection*{Reduction from $C_3$ to $A_2$.}
\begin{gather*}
{\rm Pr}=\left(\begin{array}{ccc}
1&1&2\\
0&1&0
\end{array}\right);
\\[2ex]
\begin{array}{@{}l}
 C_{(1,0,0)}(C_3)\mapsto C_{(1,0)}(A_2)+C_{(0,1)}(A_2), \\
      C_{(0,1,0)}(C_3)\mapsto C_{(1,1)}(A_2)+C_{(1,0)}(A_2)+C_{(0,1)}(A_2),\\
      C_{(0,0,1)}(C_3)\mapsto C_{(2,0)}(A_2)+C_{(0,2)}(A_2)+C_{(0,0)}(A_2);
\end{array}\qquad \\ \text{therefore} \qquad \begin{array}{l}
  X_1\mapsto\ Y_1+Y_2,\\
  X_2\mapsto\ Y_1Y_2+Y_1+Y_2-3,\\
  X_3\mapsto\ Y_1^2+Y_2^2-2Y_1-2Y_2+2.\end{array}\end{gather*}

\noindent{Example of reduction: }
In this case, $C_{(a,b,c)}(\overline{Y}):=C_{(a,b)}(Y_1,Y_2)$  are polynomials of the subalgebra $A_2$.
\begin{gather*}
  C_{(1,0,1)}(C_3) = X_1X_3-2X_2 \mapsto Y_1^3+Y_2^3+Y_1^2Y_2+Y_1Y_2^2-2Y_1^2-2Y_2^2-6Y_1Y_2+6=\\=C_{(3,0)}(\overline{Y})+C_{(0,3)}(\overline{Y})+C_{(2,1)}(\overline{Y})+C_{(1,2)}(\overline{Y})+Y_1+Y_2;\\
  C_{(2,0,0)}(C_3) = X_1^2-2X_2-6 \mapsto Y_1^2+Y_2^2-2Y_1-2Y_2=C_{(2,0)}(\overline{Y})+C_{(0,2)}(\overline{Y}).
\end{gather*}

\subsection*{Reduction from $C_3$ to $A_1$.}
\begin{gather*}
{\rm Pr}=\left(\begin{array}{ccc}5&8&9  \end{array}\right);
\\[2ex]
\begin{array}{@{}l}
C_{(1,0,0)}(C_3)\mapsto C_5(A_1)+C_3(A_1)+C_1(A_1),\\
      C_{(0,1,0)}(C_3)\mapsto C_8(A_1)+C_6(A_1)+2C_4(A_1)+2C_2(A_1),\\
      C_{(0,0,1)}(C_3)\mapsto C_9(A_1)+C_7(A_1)+C_3(A_1)+C_1(A_1);
\end{array}\qquad \\ \text{therefore} \qquad \begin{array}{l}
  X_1\mapsto\ Y^5-4Y^3+3Y,\\
  X_2\mapsto\ Y^8-7Y^6+16Y^4-13Y^2,\\
  X_3\mapsto\ Y^9-8Y^7+20Y^5-15Y^3.\end{array}\end{gather*}

\noindent{Example of reduction: }
\begin{gather*}
C_{(1,1,0)}(C_3) = X_1X_2 {-} 3 X_3 {-} 4X_1 \mapsto
Y^{13}{-}11Y^{11}{+}44Y^9{-}74Y^7{+}36Y^5{+}22Y^3{-}12Y=\\=\tilde{T}_{13}(Y){+}2\tilde{T}_{11}(Y)+\tilde{T}_{9}(Y){+}3\tilde{T}_{7}(Y){+}2\tilde{T}_5(Y){-}43\tilde{T}_{3}(Y){+}128Y;
\\
C_{(1,0,1)}(C_3) = X_1X_3{-}2X_2 \mapsto
Y^{14}{-}12Y^{12}{+}55Y^{10}{-}121Y^8{+}134Y^6{-}77Y^4{+}26Y^2=\\=\tilde{T}_{14}(Y)+2\tilde{T}_{12}(Y)+2\tilde{T}_{10}(Y)+\tilde{T}_8(Y)+2\tilde{T}_6(Y){+}\tilde{T}_4(Y)+\tilde{T}_2(Y)+2.
\end{gather*}

\section{Concluding remarks}

Given the branching rules for polynomials, their applications now can be contemplated, and conceivably uncommon applications can be stimulated.
When changing orbit functions to polynomials, the substitution wipes out explicit dependence on group parameters.
In spite of that, the structure of polynomials of a given family `knows' about its Lie group. It is uniquely determined by that group.

Perhaps the most important exploitation of Chebyshev polynomials  of one variable is the optimal interpolation of functions, and evaluation of their integrals.
We know \cite{LX} that these results extend to polynomials of groups of type $A_n$. Clearly the other types of simple Lie groups need to be investigated as well.

The relation between Chebyshev and Jacobi polynomials is known in one variable.
In multivariate context this relation can be established using the fact that our polynomials in $n$ variables are
modification of the monomial symmetric polynomials and the last ones are the building blocks of Jacobi polynomials in $n$ variables
(see Sec.~11 of~\cite{KP1}).

\section*{Acknowledgements}
JP gratefully acknowledges support of this research by the Natural Sciences and Engineering Research Council of Canada and by the MIND Research Institute of Santa Ana, California.
All authors express their gratitude for the hospitality extended to them at the Doppler Institute of the Czech Technical University, where most of the work was done.


\newpage
\section*{Appendix}\label{appendix}
There are three simple Lie groups of rank three, namely $A_3$, $B_3$ and $C_3$.
In this section we present all the numerical information necessary to deal with
orbit functions of three variables and polynomials in three variables.
A part of these data were taken from~\cite{NP2008}, where some more information about groups of ranks two and three can be found.
\medskip

{\bf Cartan matrices:}
\begin{gather*}
A_2: \
\mathfrak{C}=\left(
\begin{smallmatrix}
2&-1\\
-1&2
\end{smallmatrix}\right) ,\qquad
C_2: \
\mathfrak{C}=\left(\begin{smallmatrix}
2&-1\\
-2&2
\end{smallmatrix}\right),
\\
A_3: \
\mathfrak{C}=\left(
\begin{smallmatrix}
2&-1&0\\
-1&2&-1\\
0&-1&2
\end{smallmatrix}\right) ,\qquad
B_3: \
\mathfrak{C}=\left(\begin{smallmatrix}
2&-1&0\\
-1&2&-2\\
0&-1&2
\end{smallmatrix}\right) ,\qquad
C_3: \
\mathfrak{C}=\left(\begin{smallmatrix}
2&-1&0\\
-1&2&-1\\
0&-2&2
\end{smallmatrix}\right).
\end{gather*}
\medskip

{\bf Weyl group orders:}
\begin{gather*}
|W(A_n) |= (n+1)!,\qquad |W(B_n)|=|W(C_n)|=2^n n!.
\end{gather*}
\medskip

{\bf Number of points in the Weyl group orbit:}\\[2ex]
\begin{center}
\begin{tabular}{|c|c|c|c|c|c|c|c|}
\cline{6-8}
\multicolumn{5}{c|}{\quad}                                  &$\lambda$             &$A_3$ & $B_3$ and $C_3$\\
\cline{6-8}
\multicolumn{5}{c|}{\quad}                                  &$(\star,\star,\star)$ &24    &48\\
\cline{1-4}\cline{6-8}
$\lambda$             & $A_2$ & $C_2$ & $G_2$  &\qquad\qquad&$(\star,\star,0)$     &12    &24\\
\cline{1-4}\cline{6-8}
$(\star,\star)$       &   6   &   8   &   12               &&$(\star,0,\star)$     &12    &24\\
\cline{1-4}\cline{6-8}
$(\star,0)$           &   3   &   4   &    6               &&$(0,\star,\star)$     &12    &24\\
\cline{1-4}\cline{6-8}
$(0,\star)$           &   3   &   4   &    6               &&$(\star,0,0)$         &4     &6\\
\cline{1-4}\cline{6-8}
\multicolumn{5}{c|}{\quad}                                  &$(0,\star,0)$         &6     &12\\
\cline{6-8}
\multicolumn{5}{c|}{\quad}                                  &$(0,0,\star)$         &4     &8\\
\cline{6-8}
\end{tabular}
\end{center}
\medskip

{\bf Highest roots:}
\begin{gather*}
\xi_{A_3}=\alpha_1+\alpha_2+\alpha_3,\qquad
\xi_{B_3}=\alpha_1+2\alpha_2+2\alpha_3,\qquad
\xi_{C_3}=2\alpha_1+2\alpha_2+\alpha_3.
\end{gather*}
\medskip

{\bf Fundamental regions:}
\begin{gather*}
F(A_3)=\{0,\check{\w}_1,\check{\w}_2,\check{\w}_3\},\qquad
F(B_3)=\{0,\ \check{\w}_1,\ \tfrac12 \check{\w}_2,\ \tfrac12 \check{\w}_3\},\qquad
F(C_3)=\{0,\tfrac12\check{\w}_1,\tfrac12\check{\w}_2,\check{\w}_3\}.
\end{gather*}
\medskip

{\bf Discretizations of fundamental regions:}
\begin{gather*}
F_M(A_3)=\left\{\frac{s_1}{M}\cw_1+\frac{s_2}{M}\cw_2+\frac{s_3}{M}\cw_3\mid M,s_1,s_2,s_3\in \Z^{\ge 0}, s_1+s_2 +s_3 \le M\right\},\\
F_M(B_3)=\left\{\frac{s_1}{M}\cw_1+\frac{s_2}{M}\cw_2+\frac{s_3}{M}\cw_3\mid M,s_1,s_2,s_3\in \Z^{\ge 0}, 2s_1+2s_2 +s_3 \le M\right\},\\
F_M(C_3)=\left\{\frac{s_1}{M}\cw_1+\frac{s_2}{M}\cw_2+\frac{s_3}{M}\cw_3\mid M,s_1,s_2,s_3\in \Z^{\ge 0}, s_1+2s_2 +2s_3 \le M\right\}.
\end{gather*}
\begin{gather*}
\widetilde{F}_M(A_3)=\left\{\frac{s_1}{M}\cw_1+\frac{s_2}{M}\cw_2+\frac{s_3}{M}\cw_3\mid M,s_1,s_2,s_3\in \N, s_1+s_2 +s_3 \le M\right\},\\
\widetilde{F}_M(B_3)=\left\{\frac{s_1}{M}\cw_1+\frac{s_2}{M}\cw_2+\frac{s_3}{M}\cw_3\mid M,s_1,s_2,s_3\in \N, 2s_1+2s_2 +s_3 \le M\right\},\\
\widetilde{F}_M(C_3)=\left\{\frac{s_1}{M}\cw_1+\frac{s_2}{M}\cw_2+\frac{s_3}{M}\cw_3\mid M,s_1,s_2,s_3\in \N, s_1+2s_2 +2s_3 \le M\right\}.
\end{gather*}
\medskip

{\bf Number of points in $F_M$ and $\widetilde{F}_M$:}
\begin{gather*}
|F_M(A_n)|={   M+n \choose n};\\
|F_{2k}(B_n)|= |F_{2k}(C_n)|= {n+k \choose n}+{  n+ k-1 \choose n },\quad k\in \N; \\
|F_{2k+1}(B_n)|= |F_{2k+1}(C_n)|= 2{n+k \choose  n},\quad k\in \N.
\end{gather*}

\begin{gather*}
|\widetilde{F}_M(G)|=\left\{
\begin{array}{cll}
0&\textrm{for}&M<h,\\
1&\textrm{for}&M=h,\\
|F_{M-m}|\;&\textrm{for}&M>h,
\end{array}\right.
\qquad \text{where} \quad h \quad \text{is the Coxeter number.}
\end{gather*}
\medskip

{\bf Weight functions $f$ and $\tilde{f}$ for discrete orthogonality of $C$- and $S$-polynomials:}\\[2ex]
$f(A_2){=}|W_u(A_2)|J(u)$,
where $J(u)^{-1}{=}4\pi^2\sqrt{u_1^2u_2^2{-}4u_1^3{-}4u_2^3{+}18u_1u_2{-}27}$.
\medskip

\noindent
$f(C_2){=}|W_u(C_2)|J(u)$,
where $J(u)^{-1}{=}4\pi^2\sqrt{u_1^2u_2^2{-}4(u_1^4{+}4 u_1^2{-}u_2^3{-}8u_2^2{-}16u_2{+}6u_1^2 u_2)}$.
\medskip

\noindent
$f(A_3){=}\frac{|W_u(A_3)|}{8\pi^3\sqrt{|\mathbf{S}(u)|}}$\quad
and\quad
$\tilde{f}(A_3){=}\frac{4!\sqrt{|\mathbf{S}(u)|}}{8\pi^3}$,
where\\
$\mathbf{S}(u)=256{-}27u_1^4{+}144u_1^2u_2{-}128u_2^2{-}4u_1^2u_2^3{+}16u_2^4{-}192u_1u_3{+}18u_1^3u_2u_3{-}80u_1u_2^2u_3{-}6u_1^2u_3^2$\\
$\phantom{\mathbf{S}(u)=}{+}144u_2u_3^2{+}u_1^2u_2^2u_3^2{-}4u_2^3u_3^2{-}4u_1^3u_3^3{+}18u_1u_2u_3^3{-}27u_3^4.$
\medskip

\noindent
$f(B_3){=}\frac{|W_u(B_3)|}{8\pi^3\sqrt{|\mathbf{S}(u)|}}$\quad
and\quad
$\tilde{f}(B_3){=}\frac{2^3 3!\sqrt{|\mathbf{S}(u)|}}{8\pi^3}$,
where\\
$\mathbf{S}(u)=(16{+}4u_2{-}u_3^2)(1728{+}1728u_1{+}432u_1^2{-}32u_1^3{-}16u_1^4{+}864u_2{+}576u_1u_2{+}72u_1^2u_2{-}8u_1^3u_2$\\
$\phantom{\mathbf{S}(u)=}{+}108u_2^2{+}36u_1u_2^2{-}u_1^2u_2^2{+}4u_2^3{-}432u_3^2{-}216u_1u_3^2{+}4u_1^3u_3^2{-}108u_2u_3^2{-}18u_1u_2u_3^2{+}27u_3^4).$
\medskip

\noindent
$f(C_3){=}\frac{|W_u(C_3)|}{8\pi^3\sqrt{|\mathbf{S}(u)|}}$\quad
and\quad
$\tilde{f}(C_3){=}\frac{2^3 3!\sqrt{|\mathbf{S}(u)|}}{8\pi^3}$,
where\\
$\mathbf{S}(u)=(u_3{-}2u_2{+}4u_1{-}8)(8{+}4u_1{+}2u_2{+}u_3)(u_1^2u_2^2{-}4u_2^3{-}4u_1^3u_3{+}18u_1u_2u_3{-}27u_3^2)$.
\medskip

{\bf Coxeter numbers:}
\begin{gather*}
h=\left\{\begin{array}{ll}
n+1&\textrm{for } A_n\\
2n&\textrm{for } B_n \textrm{ and } C_n.
\end{array}\right.
\end{gather*}
\medskip

{\bf Level vectors:}\\[2ex]
\begin{center}
\begin{tabular}{|c|l|}
\hline
Lie group & level vector\\
\hline
$A_2$ &  $(2,2)$\\
\hline
$C_2$ & $(3,4)$\\
\hline
$A_3$ & $(3,4,3)$\\
\hline
$B_3$ & $(6,10,6)$\\
\hline
$C_3$ & $(5,8,9)$\\
\hline
\end{tabular}
\end{center}
\medskip

{\bf Weyl group orbits:}\\[2ex]
Below we present diagrams of the Weyl orbits of groups $A_2$, $C_2$, $A_3$, $B_3$ and $C_3$,
each diagram contains only half of all its points,
the rest of points can be obtained by means of change of sign in front of each coordinate of the diagram points.
\medskip

\hspace{-10ex}~$\begin{array}{cc}
{
\xymatrix{
&&*+[F]{(a,b)}\ar[dl]^{r_1}\ar[dr]_{r_2}&\\
&*+[F]{({-}a,a{+}b)}&&*+[F]{(a{+}b,{-}b)}}
}
&
{
\xymatrix{
&&*+[F]{(a,b)}\ar[dl]^{r_1}\ar[dr]_{r_2}&\\
&*+[F]{({-}a,a{+}b)}\ar[d]^{r_2}&&*+[F]{(a{+}2b,{-}b)}\ar[d]^{r_1}\\
&*+[F]{(a{+}2b,{-}a{-}b)}&&*+[F]{({-}a{-}2b,a{+}b)}}
}\\
{}&{}\\
W_{(a,b)}(A_2)& W_{(a,b)}(C_2)
\end{array}$

\begin{figure}[h]
\centerline{\includegraphics[scale=0.85]{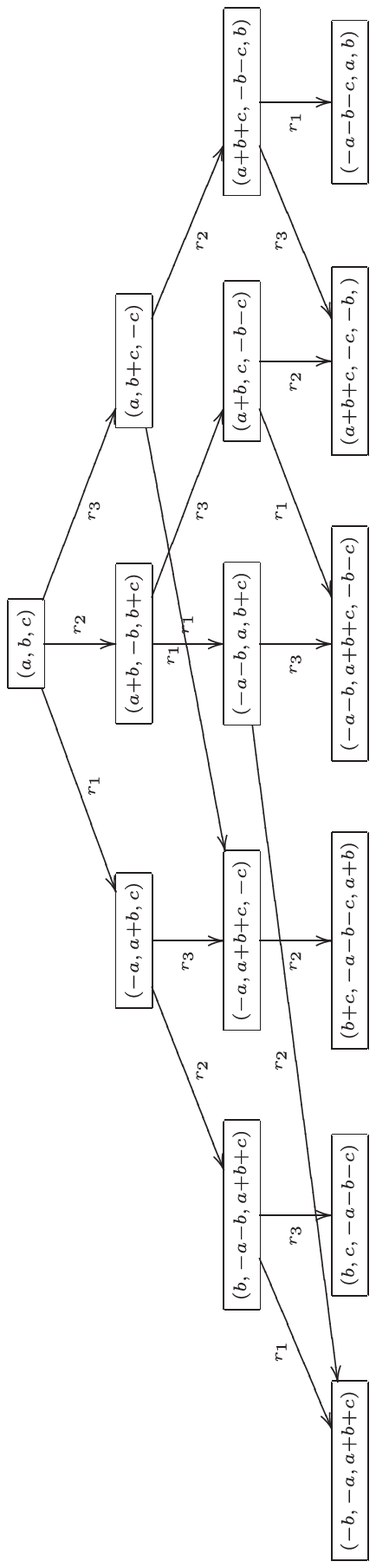}}
\caption{Weyl group orbit of a generic point $\lambda=(a,b,c)_{\w}$ in the case of group $A_3$.}\label{fig_a3orb}
\end{figure}

\begin{figure}[h]
\centerline{\includegraphics[scale=0.8]{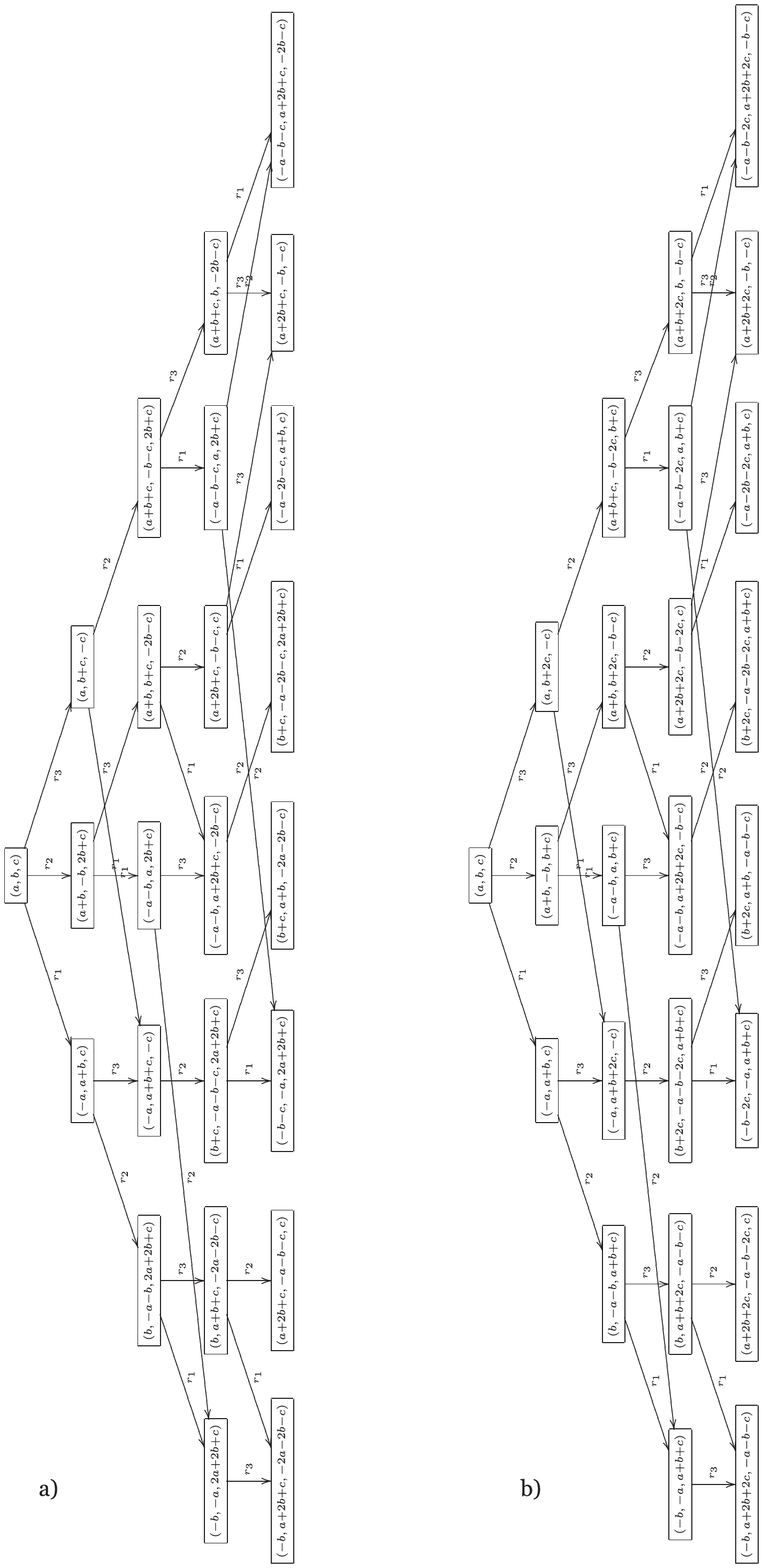}}
\caption{Weyl group orbits: a) $W_{(a,b,c)}(B_3)$\quad and\quad b) $W_{(a,b,c)}(C_3)$.}\label{fig_b3c3orb}
\end{figure}


\begin{thebibliography}{99}
\itemsep-.05ex

\bibitem{BMP}
Bremner~M.R., Moody~R.V., Patera~J.,
{\it Tables of dominant weight multiplicities for representations of simple Lie algebras},
Marcel Dekker, New York 1985, 340 pages.

\bibitem{Dunkl}
Dunkl~Ch., Xu~Yu.,
{\it Orthogonal polynomials of several variables},
Cambridge University Press, New York, 2008, 408 pages.

\bibitem{HP}
Hrivn\'ak~J., Patera~J.,
On discretization of tori of compact simple Lie Groups,
\textit{J.~Phys.~A: Math. Theor.}, \textbf{42} (2009), 385208, 26 pages, math-ph/0905.2395.

\bibitem{IlievXu}
Iliev~P., Xu~Yu.,
Discrete orthogonal polynomials and difference equations of several variables,
\textit{Adv. in Appl. Math.} \textbf{212}, (2007), 1--36.

\bibitem{KassMoodyPateraSlansky1990}
Kass~S., Moody~R.V., Patera~J., Slansky~R.,
\emph{Affine Lie algebras, weight multiplicities, and branching rules}, Vol.1 and Vol.2,
(Los Alamos Series in Basic and Applied Sciences, University of California Press, Berkeley, 1990).

\bibitem{KP1}
Klimyk~A., Patera~J.,
Orbit functions,
{SIGMA} \textbf{2} (2006), 006, 60~pages, math-ph/0601037.

\bibitem{KP2}
Klimyk~A., Patera~J.,
Antisymmetric orbit functions,
{SIGMA} \textbf{3} (2007), 023, 83~pages, math-ph/0702040v1.


\bibitem{Koor}
Koornwinder~T.H.
Orthogonal polynomials in two variables which are eigenfunctions of two algebraically independent partial differential operators I-IV,
\textit{Nedrl. Akad. Wetensch. Proc. Ser. A.}, {\bf 77}, {\bf 36} (1974) 48-66, 357-381.


\bibitem{LNP}
Larouche~M., Nesterenko~M., Patera~J.,
Branching rules for orbits of the Weyl group of the Lie algebra~$A_n$,
{\it J.~Phys,~A: Math. Theor.}, {\bf 42} (2009),  485203, 14 pages; arXiv:0909.2337.

\bibitem{LX}
Li~H., Xu~Yu.,
Discrete Fourier analysis on fundamental domain of $A_n$ lattice and on simplex in $d$-variables,
\textit{arXiv:0809.1079} (2008), 39 pages.

\bibitem{Lidl}
Lidl~R., Wells~Ch.,
Chebyshev polynomials in several variables.
\textit{Journal f\"ur die reine und angewandte Mathematik} \textbf{255} (1972), 104--111.

\bibitem{Macdonald}
Macdonald~I.G.,
Some conjectures for root systems
\textit{SIAM Journal on Mathematical Analysis} \textbf{13} (6) (1982), 988–-1007.

\bibitem{Macdonald1}
Macdonald I.G.,
Orthogonal polynomials associated with root systems,
{\it S\'eminaire Lotharingien de Combinatoire}, Article B45a, Strasbourg, 2000.

\bibitem{McP}
McKay~W.G., Patera~J.,
{\it Tables of dimensions, indices, and branching rules for representations of simple Lie algebras},
Marcel Dekker, New York, 1981, 317 pages.

\bibitem{MP87}
Moody~R.V., Patera~J.,
Computation of character decompositions of class functions on compact semisimple Lie groups,
{\it Mathematics of Computation} {\bf 48}, (1987) 799--827.

\bibitem{MP}
Moody~R.V., Patera~J.,
Orthogonality within the Families of $C$--, $S$-- and $E$-- Functions of Any Compact Semisimple Lie Group,
{SIGMA} \textbf{2} (2005), 076, 14~pages, math-ph/0512029.

\bibitem{MPi}
Moody~R.V., Pianzola~A.,
Lambda-mappings between the representation rings of Lie algebras,
\textit{Can. J. Math.}, \textbf{35} (1983), 898--960.

\bibitem{NP2008}
Nesterenko~M., Patera~J.,
Three dimensional C-, S- and E-transforms,
{\it J.Phys.A: Math. Theor.}, \textbf{41} (2008), 475205, 31 pages; arXiv:0805.3731v1.

\bibitem{NPT0}
Nesterenko~M., Patera~J., Tereszkiewicz~A.,
Orbit functions of $SU(n)$ and Chebyshev polynomials,
\emph{arXiv:0905.2925v2} (2009), 15 pages.

\bibitem{NPT}
Nesterenko~M., Patera~J., Tereszkiewicz~A.
Orthogonal Polynomials of compact simple Lie groups,
\emph{arXiv:1001.3683v2} (2010), 31 pages.


\bibitem{PSS}
Patera~J., Sharp~R.T., Slansky~R.,
On a new relation between semisimple Lie algebras,
{\it J.~Math.~Phys.}, {\bf 21}, (1980), 2335--2341.

\bibitem{P}
Patera~J.,
Compact simple Lie groups and theirs $C$-, $S$-, and $E$-transforms,
{SIGMA} \textbf{1} (2005), 025, 6~pages, math-ph/0512029.



\bibitem{R}
Rivlin~T.J.,
{\it The Chebyshef polynomials\/},
Wiley, New York, 1974.


\bibitem{Xu}
Xu~Yu.,
On discrete orthogonal polynomials of several variables,
\it{Adv. in Appl. Math.} \textbf{33}, (2004), 615--632.


\end{thebibliography}
\end{document}